\documentclass[12pt]{amsart}
\usepackage{amssymb,amsthm,amsmath,amstext}
\usepackage{mathrsfs}  % allows nice script font for sheaves
\usepackage{bm}        % allows bold italic letters in math mode
\usepackage{mathtools} % allows more extendible arrows

\usepackage[left=1.28in,top=.8in,right=1.28in,bottom=.8in]{geometry}

\usepackage{graphicx}
\usepackage[all]{xy}
\theoremstyle{plain}
\newtheorem{theorem}{Theorem}[section]
\newtheorem{prop}[theorem]{Proposition}
\newtheorem{lemma}[theorem]{Lemma}
\newtheorem{question}[theorem]{Question}
\newtheorem{cor}[theorem]{Corollary}

\theoremstyle{definition}

\newtheorem{example}[theorem]{Example}

\theoremstyle{remark}
\newtheorem{remark}[theorem]{Remark}

\newcommand{\sheaf}[1]{\mathscr{#1}}

\newcommand{\OO}{\sheaf{O}}

\newcommand{\PP}{\sheaf{P}}
\newcommand{\XX}{\sheaf{X}}
\newcommand{\YY}{\sheaf{Y}}

\newcommand{\RR}{\sheaf{R}}

 \newcommand{\rama}{\mbox{ram$_{\XX}(\alpha)$}}
\newcommand{\ramz}{\mbox{ram$_{\XX}(\zeta)$}}
\newcommand{\ramap}{\mbox{ram$_{{\XX}'}(\alpha)$}}
\newcommand{\ramzp}{\mbox{ram$_{{\XX}'}(\zeta)$}}

\newcommand{\Z}{\mathbb Z}

\renewcommand{\P}{\mathbb P}
\newcommand{\Q}{\mathbb Q}
\newcommand{\F}{\mathbb F}

\usepackage[backref=page]{hyperref}

\hypersetup{colorlinks=true,linkcolor=blue,anchorcolor=blue,citecolor=blue,letterpaper=true}

\begin{document}

\title[Degree three cohomology] 
{ Degree three cohomology of function fields of surfaces}

\author[Parimala]{R.\ Parimala }
\address{Department of Mathematics \& Computer Science \\ %
Emory University \\ %
400 Dowman Drive~NE \\ %
Atlanta, GA 30322, USA}
\email{ parimala@mathcs.emory.edu, suresh@mathcs.emory.edu}
 
\author[Suresh]{V.\ Suresh}

\date{}

\begin{abstract}
Let $F$ be the function field of a surface $X$ over a finite field.
Let $l$ be a prime not equal to the characteristic of $F$.  
Our main result is a   
local-global principle (with respect to discrete valuations 
of $F$) for a class in  $H^3(F, \mu_l^{\otimes 3})$ to be 
written as the cup-product of a given symbol in  
 $H^2(F, \mu_l^{\otimes 2})$ by some class in 
 $H^1(F, \mu_l) \simeq F^*/F^{*l}$. 
 From this  we deduce  if $F$ contains a primitive 
 $l^{\rm th}$ root of unity, then  any  element in $H^3(F, \mu_l^{\otimes 3})$ is 
 represented by a symbol,  a result which  was previously 
 known only for $l = 2$. 
The local-global principle also   enables us to show that   if a 
threefold $X$ over a finite field  $\F$ (with char$(\F) \neq 2$) 
admits a conic bundle fibration over a surface,
 then the unramified cohomology group $H^3_{nr}(\F(X)/\F, \Z/2)$ vanishes.
Combined with earlier works, this gives a proof
  of the integral Tate conjecture for one-cycles
on certain classes of threefolds over $\F$. It also leads to a proof that
the Brauer-Manin obstruction to existence of a zero-cycle of degree 1
is the only obstruction for a certain classes of  surfaces over 
global fields of positive characteristic.
\end{abstract}

\maketitle

\section{Introduction}
 
Let $k$ be a global field  or a local field.  It follows from  class field theory 
 that central division algebras over $k$ are cyclic. If $k$
contains $l^{th}$ roots of unity, for a prime $l$ not equal to the
characteristic of $k$,   every element in $H^2(k, ~\mu_l^{\otimes 2})$ is a
symbol.   A natural question is whether a higher dimensional analogue of this
result holds. More specifically,  let $F$ be a function field in one variable over $k$ ($k$
being a global field or a local field) which contains  the $l^{th}$ roots
of unity. Is every element in $H^3(F, ~\mu_l^{\otimes 3})$ a symbol? This
question, for $F = {\Q}_p(t)$ and $l = 2$ was raised by
Serre (\cite{Se}, \S 8.3) in the context of the study of $G_2$-torsors
over this field. 

Let $F$ be the function field of a $p$-adic curve.
In (\cite{PS1}, 3.9), 
we proved that every element in $H^3(F, ~\mu_2)$ is a symbol for  $p \neq 2$. 
In view of  (\cite{PS2},  \cite{hb} and \cite{L}), every quadratic
 form in at least nine variables   over the function field $F$  of  a curve  over 
a  $p$-adic field  has a
non-trivial zero.  In particular,  every element in $I^3(F)$ is represented 
by a 3-fold Pfister form  (cf. \cite{PS2},  4.1) and  by the surjectivity  
of the Arason invariant  (\cite{MS}), every element in $H^3(F,
~\mu_2^{\otimes 3})$ is a symbol.
 In ( \cite{PS2}), we proved that every
element in $H^3(F, ~\mu_l^{\otimes 3})$ is a symbol if  $l \neq p$
and $F$ contains a primitive $l^{\rm th}$ root of unity. This was a key
ingredient in the final proof (\cite{PS2}) that over a function field in
one variable over a non-dyadic $p$-adic field, every quadratic form
in at least nine variables has a non-trivial zero.

We observe that if $k$ is a global field of positive characteristic
$p \neq 2$ and $F$ a function field in one variable over $k$, $F$ is
a $C_3$ field and hence every quadratic form in 9 variables over $F$ has a
non-trivial zero. In particular every element in $H^3(F, ~\mu_2^{\otimes 3})$ is
a symbol. In this paper, we prove that if $l$  is a prime not equal
to $p$ and $F$ contains a primitive $l^{th}$ root of unity, then
every element in $H^3(F, ~\mu_l^{\otimes 3})$ is a symbol. 
The main result of this paper is  
   the following 
 local-global principle (\ref{local-global-ps2}, \ref{local-global-ff2})
 which is a key ingredient in the proof that every element in 
 $H^3(F, ~\mu_l^{\otimes 3})$ is a symbol.

\begin{theorem}
\label{int-main} Let $k$ be either a $p$-adic field or a  finite field of 
characteristic $p$. Let $X$ be a
smooth projective integral variety over $k$ and $K$ the function field of $X$. 
Suppose  that $X$ is a curve if $k$ is a $p$-adic field and $X$ is a surface 
if $k$ is a finite field.  Let $ l$ be a prime  not equal to $p$.  
 Let $\alpha \in H^2(K, \mu_l^{\otimes 2})$
be a symbol and $\zeta \in H^3(K, \mu_l^{\otimes 3})$. Suppose that for every discrete valuation 
$\nu$ of $K$ there exists $f_\nu \in K_\nu^*$ such that  $\zeta =
\alpha \cdot (f_\nu) \in H^3(K_\nu, \mu_l^{\otimes 3})$, $K_\nu$ being the completion 
of $\nu$. Then there exists $f
\in K^*$ such that $\zeta = \alpha \cdot (f)$. 
\end{theorem}

Colliot-Th\'el\`ene   raised the following question (\cite{CT3},  
cf.   \cite{CTK}, Question 5.4).

\begin{question}
Let $k$ be a  finite field and $l$ a prime not equal to the characteristic of
$k$. Let $X$ be a smooth projective 3-fold over $k$ and $K$ its function field.
Is   $H^3_{nr}(K/k, \Q_l/\Z_l(2)) = 0$?
\end{question}

The above   local-global
principle in the setting of function fields of curves over global fields of 
positive  characteristic not equal to 2  leads to an affirmative answer to the
above question  for  conic fibrations over  surfaces over finite fields (\ref{h3ur-ff}).
In fact we prove the following 

 \begin{theorem}
 \label{vanishing}
Let $k$ be a finite  field of characteristic
not equal to 2. Let $X$ be a smooth, projective, integral surface
over $k$ and $K$ its function field. Let $n$ be a natural number
coprime to the characteristic of $K$ and $C$ a smooth conic over $K$. Then
$H^3_{nr}(K(C)/k, ~\mu_n^{\otimes 2}) = 0$.
\end{theorem}

Let $X$ be a smooth  projective surface over a finite field
${ {\F}}$ of characteristic $p$ not equal to 2. Let $f : Y \to
X$ be a surjective morphism, $Y$ being a smooth, projective 3-fold
over ${ {\F}}$, with the generic fibre of $f$ a smooth conic
over ${{\F}}(X)$. The vanishing of  the third  unramified cohomology
group has consequences in the study of integral Tate conjecture and
Brauer-Manin obstruction to the  existence of zero-cycles of degree
one  (cf. \cite{CTK}, Theorem 7.4). 

The theorem (\ref{vanishing}) implies that  if $X$ is a smooth projective
geometrically  integral ruled surface over ${\F}$
 and $Y \to X$ a surjective morphism with the generic fibre a
smooth conic and $Y$ smooth projective over $F$, 
then the cycle map $CH^2(Y) \otimes {{\Z}}_l \to
H^4(Y, {{\Z}}_l(2))$  is surjective  (\ref{cyclemap}). If $X = C \times
{{\P}}^1$, $C$ being a smooth projective curve and $Y \to C$
the composite morphism $Y \to C \times { {\P}}^1 \to C$, then
the Brauer-Manin obstruction is the only obstruction to local-global
principle for the existence of zero-cycles of degree one on the
generic fibre $Y_{\eta}$, $\eta$ being  the generic point of $C$
(\ref{bmob}).

The main idea of the proof of  Theorem (\ref{int-main}) is 
the following: We choose a regular model $\XX$ of $K$ proper
over $R$, where
$R$ is either a finite field or the ring of integers in a $p$-adic field. 
We assume, after possibly replacing $\XX$ by a   blow-up,  that the ramification loci 
  $\RR$
of $\alpha$ and $\zeta$ are a union of regular curves with normal crossings. 
Given $f_x \in K_x^*$ for each codimension one point in $\RR$ 
 such that $\zeta = \alpha \cdot (f_x) $
in $H^3(K_x, \mu_l^{\otimes 3})$, we choose $f \in K^*$ such that $f = f_x$ modulo $K_x^{*l}$
for all codimension one points in $\RR$.  Then $\zeta - \alpha \cdot (f)$ is unramified 
at  all codimension one points in $\RR$. 
 Let $y \in \XX$ be a codimension one point which is
not  $\RR$.  Then,   the residue of $ \zeta - \alpha \cdot (f)$ at $y$ is 
$ \overline{\alpha}^{\nu_y(f)}$,
where $\nu_y$ is the valuation at $y$ and $\overline{\alpha}$ is the specialization of $\alpha$
at $y$.  We need to  modify $f$ such that either $\nu_y(f)$ is a multiple of $l$
or $\overline{\alpha} = 0$. 
The exhaustive study of ramification pattern of  Brauer classes on
the function field of a regular integral two-dimensional scheme, 
given by  Saltman (\cite{S2}, \cite{S3}) is used in a delicate way to modify
this  choice of $f$, leading to a proof of  the theorem.

 For   the
function field of a curve over a $p$-adic field $k$ with $p \neq l$, 
a  local-global principle, similar to Theorem (1.1),  is proved in 
(\cite{PS2}, 3.4).  An additional hypothesis is 
missed   in (\cite{PS2},  3.4);  we bridge it  in  (\ref{local-global-ps1}).  The   proof of 
 (\cite{PS2}, 4.5) uses (\cite{PS2}, 3.4) and we complete   it
  in the Appendix using  (\ref{local-global-special}).

We would like to thank Colliot-Th\'el\`ene for bringing to our
attention  questions related to the unramified cohomology of
varieties over finite fields and their connections to the
integral Tate conjecture. We also thank him for helpful
discussions when this work was carried out. We  thank
Saltman for discussions on splitting ramification on surfaces.

The first author is partially supported by National Science Foundation grant 
DMS-1001872 and the second  author is partially supported by 
National Science Foundation grant  DMS-1301785.

\section{ Some Preliminaries}

In this section we recall a few basic facts from the theory of
Galois cohomology. We refer the reader to (\cite{CT1}).

Let $F$ be a field and $l$ a prime not equal to the characteristic
of $F$. Let $\mu_l$ be the group of $l^{\rm th}$ roots of unity. For
$i \geq 1$, let $\mu_l^{\otimes i}$ be the Galois module given by
the tensor product of $i$ copies of $\mu_l$. For $n \geq 0$, let
$H^n(F, ~\mu_l^{\otimes i})$ be the $n^{\rm th}$ Galois cohomology
group with coefficients in $\mu_l^{\otimes i}$.

We have the Kummer isomorphism $F^*/F^{*l} \simeq H^1(F, ~\mu_l)$.
For $a \in F^*$, its class in $H^1(F, ~\mu_l)$ is denoted by $(a)$.
If $a_1, \cdots , a_n \in F^*$, the cup product $(a_1) \cdots (a_n)
\in H^n(F, ~\mu_l^{\otimes n})$ is  called a {\it symbol}. We have an
isomorphism $H^2(F, ~\mu_l)$ with the $l$-torsion subgroup $_lBr(F)$
of the Brauer group of $F$. We define the {\it index} of an element
$\alpha \in H^2(F, ~\mu_l)$ to be the index of the corresponding
central  division algebra in $_lBr(F)$.

For a (rank one) discrete valuation $\nu$ of  $F$, let $\kappa(\nu)$ denote the
residue field at  $\nu$, ${\OO}_{\nu}$  the ring of integers in
$F$ at $\nu$ and $F_{\nu}$ the completion of $F$ at $\nu$.  Let
$\alpha \in H^n(F, ~\mu_l^{\otimes m})$. We say that $\alpha$ is
{\it unramified} at $\nu$ if $\alpha$ is in the image of the
restriction map $H^n_{\acute{e}t}({\OO}_{\nu}, ~\mu_l^{\otimes
m}) \to H^n(F, ~\mu_l^{\otimes m})$. We say that $\alpha$ is
{\it ramified} at $\nu$ if it is not unramified at $\nu$. Suppose
char$(\kappa(\nu)) \neq l$. Then there is a {\it residue}
homomorphism 
$$ \partial_{\nu} : H^n(F, ~\mu_l^{\otimes m}) \to
H^{n-1}(\kappa({\nu}), ~\mu_l^{\otimes (m-1)}).$$ 
For  $\alpha \in
H^n(F, ~\mu_l^{\otimes m})$, $\alpha$ is unramified at $\nu$ if and
only if $\partial_{\nu}(\alpha) = 0$. Suppose  that $\alpha$ is unramified
at $\nu$.   Then, by  definition,  there exists a unique 
$\alpha_0 \in H^n_{\acute{e}t}(\OO_\nu, ~\mu_l^{\otimes m})$ which maps to $\alpha$.
 Let $\overline{\alpha} \in H^n(\kappa(\nu), ~\mu_l^{\otimes m})$ be the image of 
 $\alpha_0$ under the natural map $H^n_{\acute{e}t}(\OO_\nu, ~\mu_l^{\otimes m})
 \to H^n(\kappa(\nu), ~\mu_l^{\otimes m})$.  We call $\overline{\alpha}$
 the {\it specialization} of $\alpha$
at $\nu$. 
Let $\pi \in F^*$ be a parameter at $\nu$ and $\zeta =
\alpha \cdot (\pi) \in H^{n+1}(F, ~\mu_l^{\otimes (m+1)})$. Then 
$\overline{\alpha} = \partial_{\nu}(\zeta) \in H^n(\kappa(\nu), ~\mu_l^{\otimes m})$.

Let $\XX$ be a regular integral scheme of dimension $d$, with
function field  $F$.  For each point $x \in {\XX}$, let $A_x$ be
the local ring at $x$, $\kappa(x)$ the residue field at $x$,
$\hat{A}_x$ the completion of $A_x$ at its maximal ideal $m_x$ and
$K_x$ the field of fractions of $\hat{A}_x$.  Suppose that $l$ is a unit on ${\XX}$. Let
$\XX^1$ be the set of points of $\XX$ of codimension 1. A point
$x \in \XX^1$  gives rise to a discrete valuation $\nu_x$ on $F$.
The residue field of this discrete valuation ring is denoted by
$\kappa(x)$. The corresponding residue homomorphism is denoted by
$\partial_x$. Let $\zeta \in H^n(F,
~\mu_l^{\otimes m})$ be any element. 

We say that 

\noindent
 $\bullet$  $\zeta$ is {\it unramified} at $x$ if  $\partial_x(\zeta) = 0$;

\noindent
 $\bullet$  $\zeta$ is {\it ramified} at $x$ if  $\partial_x(\zeta) \neq 0$;

We note that   $\zeta$ is ramified only at finitely many codimension  
one points of ${\XX}$.  We define the {\it ramification divisor} $\mbox{
ram}_{\XX}(\zeta)$ to be  $\sum x$ as $x$ runs over the points of
$\XX^1$ where $\zeta$ is ramified.  The unramified cohomology on
$\XX$, denoted by $H^n_{nr}(F/\XX, ~\mu_l^{\otimes m})$, is
defined by
$$H^n_{nr}(F/\XX, ~\mu_l^{\otimes m}) = 
 \cap_x ker(\partial_x : H^n(F, ~\mu_l^{\otimes m}) \to H^{n-1}(\kappa(x),
~\mu_l^{\otimes (m-1)})),$$
  where $x$'s run  over  $\XX^1$. 
  If $\XX$ is a projective variety over a field $k$, we also
  denote $H^n_{nr}(F/\XX, ~\mu_l^{\otimes m})$ by
  $H^n_{nr}(F/k, ~\mu_l^{\otimes m})$.
 We say that
$\zeta \in H^n(F, ~\mu_l^{\otimes m})$ is {\it unramified on
$\XX$} if $\zeta \in H^n_{nr}(F/\XX, ~\mu_l^{\otimes m})$; 
if
$\XX = {\rm Spec}(R)$, then we  say that $\zeta$ is {\it unramified on}
$R$ if it is unramified on $\XX$. If $R$ is a local ring at some
point $P$ of $\XX$, then we say that $\alpha$ is {\it unramified at} $P$
if it is unramified on $R$.   Suppose $C$ is an irreducible
subscheme of ${\XX}$ of codimension 1. Then the generic point $x$
of $C$ belongs to $\XX^1$ and we set $\partial_C$ = $\partial_x$.
If $\alpha \in H^n(F, ~\mu_l^{\otimes m})$ is unramified at $x$,
then we say that $\alpha$ is {\it unramified} at $C$. The group of
elements of $H^n(F, ~\mu_l^{\otimes m})$ which are unramified at all
discrete valuations of $F$ is denoted by $H^n_{nr}(F, ~\mu_l^{\otimes
m})$.

Let $\XX$ be a  regular integral surface.
Let $D$ be a divisor on ${\XX}$. Then by the resolution of
singularities for surfaces (cf. \cite{Li1} and \cite{Li2}), there exists a
regular, integral surface  $\XX'$ with a proper birational
morphism $\XX' \to \XX$, such that the total transform of $D$ is a
union of regular curves with normal crossings (cf. \cite{Sh}, Theorem,
p.38 and Remark 2, p. 43). We use this result throughout this paper
without further reference.

 Let $K$ be a field and $\ell$ a prime not equal to the characteristic 
of $K$. Suppose that $K$ contains a primitive $\ell^{\rm th}$ root of unity. 
We fix a generator $\rho$ for the cyclic group $\mu_l$ and identify the
Galois modules $\mu_l^{\otimes i}$ with $\mu_l$. This leads to an
identification of $H^n(F, ~\mu_l^{\otimes m})$ with $H^n(F,
~\mu_l)$. The element in $H^n(F, ~\mu_l)$ corresponding to the
symbol $(a_1) \cdots (a_n) \in H^n(F, ~\mu_l^{\otimes n})$ through
this identification is again denoted by $(a_1) \cdots (a_n)$.
If $\nu$ is a discrete valuation on $K$ and $a, b \in K^*$, we have
$\partial_\nu((a) \cdot (b)) =  (\overline{\frac{a^{\nu(b)}}{b^{\nu(a)}}})$. 
If $u \in K^*$ is a unit at $\nu$, then $\partial_{\nu}((u) \cdot (a) \cdot (b))
= (\overline{u}) \cdot \partial_\nu((a) \cdot (b))$.

We now recall some facts concerning ramifications of division
algebras on surfaces from (\cite{S1}, \cite{S2}, \cite{S3}). Let $K$ be the function
field of a regular integral surface $\XX$.  Let  $C_1, \cdots , C_n$ be 
distinct irreducible curves on $\XX$ and $P \in \XX$ be a closed point.

\noindent
$\bullet$ We say that $P$ is {\it curve point} of the curves 
${C_1, \cdots , C_n}$ 
if $P \in C_i$
for some $i$ and $P \not\in C_j$ for all $j \neq i$. 

\noindent
 $\bullet$ We say that $P$ is a {\it nodal point} of the curves 
${C_1, \cdots , C_n}$ if $P \in C_i \cap C_j$ for some $i \neq j$.

 Let $l$ be a prime which is a unit on ${\XX}$. 
Suppose that $K$ contains a primitive $l^{\rm th}$ root  of unity. 
Let $\alpha \in H^2(K, ~\mu_l)$ be such that
ram$_{\XX}(\alpha)$ is a union of regular curves with normal
crossings.  Let $P \in {\XX}$ be a closed point. 

\noindent
 $\bullet$ If $P$ is a curve point
of the support of ram$_{\XX}(\alpha)$, we
say that $P$ is a {\it curve point} for $\alpha$. 

\noindent
 $\bullet$ If $P$ is a nodal point of the support of ram$_{\XX}(\alpha)$,
then we say that $P$ is a {\it nodal point} for $\alpha$.

 If $P$ is not on any curve in the support of ram$_{\XX}(\alpha)$,  
then $\alpha$ is unramified on the local ring ${\OO}_{\XX, P}$. 
Suppose that $P$ is a curve point for $\alpha$. 
Let $C$ be the irreducible curve in  the support of ram$_{\XX}(\alpha)$ and
   $\pi \in {\OO}_{\XX, P}$ define the curve $C$ at $P$.  By 
(\cite{S1}, 1.1) $\alpha = \alpha' + (u) \cdot
(\pi)$ where $\alpha' \in H^2(K, ~\mu_l)$ is unramified at $P$ and $u$ is 
a unit at $P$.  In particular,  $\partial_C(\alpha) = (\overline{u}) \in
H^1(\kappa(C), ~\mu_l)$.   Further,  if $w$ is a unit at $P$ with
$\partial_C(\alpha) = (\overline{w})  \in 
H^1(\kappa(C), ~\mu_l)$, then $\alpha - (w) \cdot (\pi)$ is unramified
at $P$.

 Suppose that $P$ is a nodal point of  
ram$_{\XX}(\alpha)$. Since the support of ram$_{\XX}(\alpha)$ is
a union of regular curves with normal crossings, $P$ is only on two
curves, say $C$ and $E$, in the support of ram$_{\XX}(\alpha)$.
 Let $\pi$ and $\delta$ in ${\OO}_{\XX, P}$  define $C$ and $E$ at $P$
respectively so that  the maximal ideal of ${\OO}_{\XX, P}$ is
generated by $\pi$ and $\delta$.  Suppose that $\partial_C(\alpha)
\in H^1(\kappa(C),~\mu_l)$ and $\partial_E(\alpha) \in H^1(\kappa(E),
~\mu_l)$ are unramified at $P$. Let $u(P), v(P) \in H^1(\kappa(P),
~\mu_l)$ be the specializations at $P$ of $\partial_C(\alpha)$ and
$\partial_E(\alpha)$ respectively. The condition that
$\partial_C(\alpha) \in H^1(\kappa(C), \mu_l)$
and $\partial_E(\alpha) \in H^1(\kappa(E),
\mu_l)$ are unramified at $P$ is equivalent to
the condition $\alpha = \alpha'+ (u) \cdot (\pi) + (v) \cdot
(\delta)$ for some units $u, v \in {\OO}_{{\XX}, P}$ and
$\alpha'$ unramified on ${\OO}_{\XX, P}$  (\cite{S2}, \S 2). The
specialisations of $\partial_C(\alpha)$ and $\partial_E(\alpha)$ in
$H^1(\kappa(P), \mu_l) \simeq
\kappa(P)^*/\kappa(P)^{*^l}$ are given by the images of $u$ and $v$
in $\kappa(P)$. 

 Following Saltman (\cite{S2}, \S 2), we say that 

\noindent
$\bullet$
 $P$ is a {\it cool point} if $u(P), v(P)$ are zero in
$H^1(\kappa(P), \mu_l)$,

\noindent
 $\bullet$ a {\it chilly point} if $u(P)$ and $v(P)$ 
generate the same subgroup of $H^1(\kappa(P), \mu_l)$ and neither
of them is zero,

\noindent
$\bullet$ a {\it hot point}  if $u(P)$ and $v(P)$ do 
not generate the same subgroup of $ H^1(\kappa(P), \mu_l)$.

\noindent
  $\bullet$ If $\partial_C(\alpha)$ and $\partial_E(\alpha)$
 are ramified at $P$, then we say that  $P$ is a {\it cold point}.

Saltman (\cite{S2}, p. 832) defined a graph using chilly
points and designated  the loops in this graph as {\it  chilly  loops}. 
He  shows  that we can blow up ${\XX}$ and assume that
there are no chilly  loops and   no cool points on $\XX$ (\cite{S2}, Corollary 2.9). 

We record here the following which we will be using frequently. 

\begin{lemma}
\label{complex} Let $A$ be a regular local ring of dimension 2
with residue field $\kappa$ and field of fractions $K$. Suppose that
$l$ is a prime not equal to char$(\kappa)$ and $K$ contains a
primitive $l^{\rm th}$ root of unity.  Let  $\zeta \in H^3(K,
\mu_l)$. Suppose that $\zeta$  is unramified on $A$ except at $(\pi)$
and $(\delta)$ for some primes $\pi$ and $\delta$ in $A$ with the maximal
ideal of $A$ generated by $\pi$ and $\delta$.  Then $\partial_{\overline
{\delta}}( \partial_{(\pi)}(\zeta) )
+ \partial_{\overline{\pi}}( \partial_{(\delta)}(\zeta))    = 0  \in  H^1(\kappa,
\mu_l)$.
\end{lemma}

\begin{proof}  By (\cite{Ka}, 1.7), the residue homomorphisms give
a complex
$$  H^3(K, \mu_l) \to \oplus_{x \in A^1}H^2(\kappa(x), \mu_l) \to
H^1(\kappa, \mu_l).$$
Since $\zeta$ is ramified possibly on $A$ only at $(\pi)$ and
$(\delta)$, we have 
$$
 \partial_{\overline{\delta}}( \partial_{(\pi)}(\zeta) )
+ \partial_{\overline{\pi}}( \partial_{(\delta)}(\zeta)) = 0  \in H^1(\kappa,
\mu_l).
$$  
\end{proof}

\begin{lemma} 
\label{purity}
Let $A$ be a regular local ring of dimension 2
with residue field $\kappa$ and field of fractions $K$. Suppose that
$l$ is a prime not equal to char$(\kappa)$ and $K$ contains a
primitive $l^{\rm th}$ root of unity.   Let $\XX$ be the simple blow-up of 
Spec$(A)$ at its maximal ideal.    \\
1) If  $\alpha  \in H^2(K, \mu_l)$ is unramified on $A$, then $\alpha$ is
unramified  on $\XX$. \\
2) If $\zeta \in H^3(K, \mu_l)$ is unramified
on $A$ and $H^2(\kappa, \mu_l) = 0$, then $\zeta$ is unramified  
on $\XX$. 
\end{lemma}

\begin{proof} 
 Suppose that $\alpha \in H^2(K, \mu_l)$ is 
unramified on $A$. 
Since $A$ is a regular ring of dimension 2, by the purity for the
Brauer group, $\alpha$ is in the image of $H^2_{\acute{e}t}(A, \mu_l)$.
 Since for any codimension one point $x$ of $\XX$,   
the local ring $\OO_{\XX, x}$ dominates $A$, $\alpha$
is in the image of $H^2_{\acute{e}t}(\OO_{\XX, x}, \mu_l)$ and hence
$\alpha$ is unramified at $x$.

 Suppose that $H^2(\kappa, \mu_l) = 0$ and $\zeta \in H^3(K, \mu_l)$ is
unramified on $A$.   Let $x \in \XX$ be a codimension one point.
Suppose that the image $y$ of $x$ in Spec$(A)$ is  a codimension one point.
Then $\OO_{\XX, x} \simeq A_y$. Since $\zeta$ is unramified at $y$, it is unramified at 
$x$. Suppose that the image of $x$ is the closed point of Spec$(A)$.
Then $x$ is the generic point of the exceptional curve  $E$ in $\XX$.
 By (\cite{Ka}, 1.7), the residue homomorphisms give
a complex
 $$  H^3(K, \mu_l) \to \oplus_{x \in \XX^1}H^2(\kappa(x), \mu_l) \to \oplus_{y\in \XX^2}
H^1(\kappa(y), \mu_l).$$
 Since $\XX^1 = $Spec$(A)^1 \cup {E}$ and $\zeta$ is unramified 
on $A$, $\partial_{{E}}(\zeta)$ is unramified at every closed point of ${E}$.
Since ${E}$, being the exceptional curve, is the projective line over $\kappa$, 
$\partial_{{E}}(\zeta)$ is in the image of $H^2(\kappa, \mu_l) \to 
H^2(\kappa({E}), \mu_l)$.
Since $H^2(\kappa, \mu_l) = 0$, $\partial_{{E}}(\zeta) = 0$ and hence $\zeta$
is unramified at ${E}$ on $\XX$.
\end{proof}

 \begin{remark}
In view of the purity theorem for higher cohomology groups due to
Gabber (cf. \cite{gabber},  Chapitre XVI, Thm. 3.1.1), 
one can dispense with the the condition $H^2(\kappa, \mu_l) = 0$ in the
second part of the above lemma. 
\end{remark}

\begin{lemma}
\label{ramified-at-pi} Let $A$ be a complete regular local ring of
dimension 2 with field of fractions $K$ and   residue field
$\kappa$.  Suppose that $l$ is a prime not equal to char$(\kappa)$ and
$K$ contains a primitive $l^{\rm th}$ root of unity. 
Let $\zeta \in H^3(K, \mu_l)$.  Suppose that $\zeta$ is unramified on
$A$ except possibly at a regular prime $(\pi)$ of $A$. If $H^2(\kappa, \mu_l) = 0$, 
then $\zeta$ is unramified on $A$. 
\end{lemma}

\begin{proof}  Let $\pi \in A$ be a regular prime. Suppose that
$\zeta$ is unramified except at $(\pi)$. By (\ref{complex}),
$\partial_{\pi}(\zeta) \in H^2(\kappa(\pi), \mu_l)$ is unramified at the maximal ideal of
$A/(\pi)$. Since $A/(\pi)$ is a complete discrete valuation ring with
 residue field $\kappa$ and $H^2(\kappa, \mu_l) =  0$,   $\partial_{\pi} (\zeta) = 0$.  
 Hence $\zeta$ is unramified on $A$. \end{proof}

For a regular integral scheme $\XX$ with function field $K$ and $P \in \XX$,
let $A_P$ denote the local ring at $P$,  $m_P$ the maximal ideal
of $A_P$, $\hat{A}_P$ the completion of $A_P$ at $m_P$ and $K_P$ the 
field of fractions of $\hat{A}_P$. 

\begin{lemma}
\label{saltman-2dim}
 Let ${\XX}$ be a regular integral surface and $K$ its 
function field.  Let $l$ be a prime which is a unit on ${\XX}$. Suppose that 
$K$ contains a primitive  $l^{\rm th}$ root of unity and $H^2(\kappa(P), \mu_l) = 0$ for every
closed point $P$ of ${\XX}$.  Let $\alpha \in H^2(K, \mu_l)$ 
be such that $\rama$ is a union of regular curves with normal crossings. Let $P \in {\XX}$ be a 
nodal point of $\alpha$ and $m_P = (\pi, \delta)$ with $\alpha$ ramified only along $\pi$ and  
$\delta$ on $A_P$.

\noindent
a) If $P$  is a chilly point of $\alpha$, then $\alpha = (u) \cdot (\pi\delta^s) \in H^2(K_P, \mu_l)$  
for  some unit $u \in A_P$ and $1 \leq s \leq l-1$.

\noindent
b)  If $P$ is a cold point of $\alpha$, then $\alpha = (u\delta) \cdot (v\pi^s) \in H^2(K_P, \mu_l)$ 
for some units $u, v \in A_P$ and $1 \leq s \leq l-1$. 

\noindent
c) If $P$ is a cool point of $\alpha$, then $\alpha = 0 \in H^2(K_P, \mu_l)$.
\end{lemma}

\begin{proof}  Suppose $P$ is a chilly point of $\alpha$. Then
$\alpha = \alpha' + (u) \cdot (\pi) + (v) \cdot (\delta)$ for some
units $u, v \in A_P$ and $\alpha' \in H^2(K, \mu_l)$ unramified at $P$
(\cite{S2}, 2.1).  Further $v(P) = u(P)^s$ for some $s$ with $1 \leq s \leq
l-1$ and hence $v = u^s \in \hat{A} _P^*/\hat{A}_P^{*^l}$.  Thus $ (v)
\cdot (\delta) = (u^s) \cdot (\delta) = (u) \cdot (\delta^s) \in
H^2(K_P, \mu_l)$.  Since $H^2(\kappa(P), \mu_l) = 0$, $\alpha' $ is
zero in $ H^2(K_P, \mu_l)$ and $\alpha = (u) \cdot (\pi\delta^s)$.

Suppose that $P$ is a cold point of $\alpha$. Then  $\alpha = \alpha'
+ (u\delta) \cdot (v\pi^s)$  for some units 
$u, v \in A_P$,  $ 1\leq s \leq l-1$
and $\alpha'$ unramified at $P$ (\cite{S2}, 2.1).  As above   
$\alpha'  $ is zero in $H^2(K_P, \mu_l)$
and   $\alpha = (u\delta) \cdot (v\pi^s)$. 

Suppose that $P$ is a cool point. Then $\alpha = \alpha' + (u) \cdot
(\pi) +  (v) \cdot (\delta)$ for some $\alpha'$ unramified at $P$
and   units $u,v \in A_P$  with
$u(P), v(P) \in \kappa(P)^{*^l}$ (\cite{S2}, 2.1).  Then 
$\alpha' $ is zero in $H^2(K_P, \mu_l)$   and $u, v \in \hat{A}_P^{*^l}$ so that 
$\alpha = 0 \in H^2(K_P, \mu_l)$.    \end{proof}

 The following lemmas are well known and we include them here for the
sake of completeness.  

 \begin{lemma}
\label{dvr-symbol}
Let $K$ be a field with a discrete valuation $\nu$, $\pi \in K^*$
be a parameter at $\nu$ and $l$ a prime not equal to the characteristic 
of  the residue field $\kappa$.  Let $\alpha  \in
H^2(K, \mu_l^{\otimes 2})$ be a symbol. Then $\alpha = 
(u) \cdot (v\pi^\epsilon)$ for some $u, v \in K^*$ units at $\nu$ and $\epsilon
= 0$ or 1.  
\end{lemma}

 \begin{proof}  Since $\alpha$ is a symbol, $\alpha = 
 (a) \cdot (b)$ for some $a, b \in K^*$.
 Write $a = u_0\pi^i$ and $b = v_0 \pi^j$ with $u_0, v_0 \in K^*$
units at $\nu$ and $i, j \in \Z$. Since $(a) \cdot (b) = (ac^l) \cdot (bd^l)$ for any
$c, d \in K^*$,  without loss of generality, we assume that $0 \leq i, j \leq l-1$.
Suppose $i = j =0$, then we have $\alpha = (u_0) \cdot (v_0)$
with $u_0$ and $v_0$ are units at $\nu$.  Suppose that 
$i$ or $j$ is non-zero. Since $(a) \cdot (b) = (b) \cdot (a^{-1})$, 
without loss of generality, we assume that  $1 \leq j \leq l-1$. 
Since $j$ is coprime to $l$,  there exists $i'$, $0 \leq i' \leq l-1$ such that
$j i' = -i$ modulo $l$.    Since $((-c)^{i'}) \cdot (c) =  i'(-c) \cdot (c) =  0$ for any 
$c \in K^*$, we have 
$$\alpha = (u_0\pi^i) \cdot (v_0\pi^j) = (u_0\pi^i (-v_0\pi^j)^{i'}) \cdot (v_0\pi^j)
= (u_0(-v_0)^{i'} \pi^{i + ji'}) \cdot (v_0\pi^j).$$
Since $i + ji'$ is a multiple of $l$, we have $\alpha = (u_0(-v_0)^{i'}) \cdot (v_0\pi^j)$.
Let $u' = u_0(-v_0)^{i'}$. Since $j$ is coprime to $l$, we have $v_0 = v^j$ modulo $K^{*l}$ 
for some $v \in K^*$. 
We have
$$ \alpha = (u') \cdot (v_0 \pi^j) = (u') \cdot ((v\pi)^j) =  (u'^{j}) \cdot (v\pi).$$
\end{proof}

 \begin{lemma}
\label{residue-unramified-at-P} Let $A$ be a two dimensional regular local ring
with field of fractions $K$, residue field $\kappa$ and 
maximal ideal $m = (\pi, \delta)$.  
 Let $l$ be a prime which is a unit in $A$. Let $\hat{A}$ be the completion
of $A$ at its maximal ideal and $\hat{K}$ the field of fractions of $\hat{A}$.
Let $n \geq 3$ and  $\zeta \in H^n(K, \mu_l^{\otimes d})$. If the image of $\zeta$
in $H^n(\hat{K}, \mu_l^{\otimes d})$ is unramified on $\hat{A}$, then 
$\partial_\pi(\zeta) \in H^{n-1}(\kappa(\pi), \mu_l^{\otimes d-1})$ is unramified 
at the discrete valuation of $\kappa(\pi)$ corresponding to the image $\overline{\delta}$
of $\delta$ in $\kappa(\pi)$. 
\end{lemma}

\begin{proof}    Let $\hat{\kappa}(\pi)$ be the
completion of $\kappa(\pi)$ at the discrete valuation corresponding to
$\overline{\delta}$. Then the field
of fractions of $\hat{A}/(\pi)$ is $\hat{\kappa}(\pi)$.
The lemma  follows from the following commutative diagram 
$$
\begin{array}{ccccccc}
H^n(K, \mu_l^{\otimes d}) &  \buildrel{\partial_\pi}\over{\to} & 
 H^{n-1}(\kappa(\pi), \mu_l^{\otimes d-1}) &  \buildrel{\partial_{\overline{\delta}}}\over{\to} &
  H^{n-2}(\kappa, \mu_l^{\otimes d-2}) \\
  \downarrow & & \downarrow & & \downarrow \\
  H^n(\hat{K}, \mu_l^{\otimes d}) &  \buildrel{\partial_\pi}\over{\to} & 
 H^{n-1}(\hat{\kappa}(\pi), \mu_l^{\otimes d-1}) &  \buildrel{\partial_{\overline{\delta}}}\over{\to} &
  H^{n-2}(\kappa, \mu_l^{\otimes d-2}),
\end{array}
$$ 
where the first two vertical arrows are  the restriction homomorphisms
 and the last vertical arrow is the identity map. 
\end{proof}

\section{Divisors on Surfaces}

 We fix the following notation  for this section:

\noindent
$\bullet$  ${\XX}$  an excellent  regular integral surface which is 
quasi-projective over an affine 
scheme. 

\noindent
 $\bullet$ $K$ the  function field  of ${\XX}$. 

\noindent
 $\bullet$  $l$  a prime which is a unit on ${\XX}$.

\noindent
 $\bullet$ Assume that $K$ contains a primitive $l^{th}$ root of unity.

\noindent
 $\bullet$
$ C_1, \cdots , C_n $ regular curves on ${\XX}$ with
normal crossings.

\noindent
 $\bullet$ ${\PP}_0$   a finite set of closed points
containing all the  nodal points of  $C_1, \cdots C_n$
   and at least one closed point  from each $C_i$. 

\noindent
 $\bullet$
 $E_1, \cdots, E_r$ irreducible curves on  ${\XX}$ 
 which do not pass through any point in ${\PP}_0$
 and  $E  = \sum t_jE_j$ with $ 1 \leq t_j \leq l-1$.

\noindent
 $\bullet$ ${\PP}$  a finite  set of closed points containing ${\PP}_0$,  
all the points of $C_i \cap E_j$
for all $i, j$ and at least one point from each $E_j$. 

\noindent
 $\bullet$ $B = {\OO
  }_{{\XX}, {\PP}}$ the regular semi-local ring at ${\PP}$
  and $\pi_i \in B$  a prime
defining $C_i$ on $B$ for $1 \leq i \leq n$.

For any two curves $D$ and $E$ on ${\XX}$ and $P \in D \cap E$, let
$(D \cdot E)_P$ denote the intersection multiplicity of $D$ and $E$ at
$P$.  For any 
$f \in K^*$, $(f)_{\XX}$ denotes the divisor of $f$ on ${\XX}$.

We recall the following from  (\cite{S3}).

\begin{lemma}  
\label{saltman1}
(\cite{S3}, 7.9) Suppose that $(C_i\cdot E)_P$ is a multiple
of $l$ for all $i$ and all $P \in {\PP}$.  There is a $z \in K^*$
such that $z$ is a unit  at the generic point of $C_i$, $z$ maps to an
$l^{\rm th}$ power  in $\kappa(C_i)^*$ for all $i$, $(z)_{\XX} = E
+ Z $ where $Z$ is a  divisor on ${\XX}$ and the
support of $Z$ does not contain  any  point in ${\PP}$.   
\end{lemma}

The following is a reformulation of (\cite{S3},  7.10 and  7.11).

\begin{lemma}
\label{saltman2}
 For each $i$,
$1 \leq i \leq n$, suppose   $u_i \in B/(\pi_i)$ are non-zero elements
such that $u_i(P) =
u_j(P)$ for all $P \in C_i \cap C_j$ and $u_i$ is a unit at all $P \in
{\PP} \setminus \cup _{i \neq j}C_i \cap C_j$. Then there exists $u  \in
K^*$  such that  $(u)_{\XX} = E + E'$ where $E'$ is a  
divisor on ${\XX}$ with the support of $E'$ 
not containing any   point in ${\PP}  \setminus \cup _{i \neq j}C_i
\cap C_j$, no $C_i$  or $E_j$ is in the support of $E'$  and the image of $u$ is
equal to $u_i$ in  $\kappa(C_i)^*/\kappa(C_i)^{*^l}$ for all $i$.  
\end{lemma}

\begin{proof} By (\cite{S3},  0.3(b)), there exists $v \in B$ such that $v = u_i$ 
modulo $(\pi_i)$.  
 By the choice,   $v$ is a unit at all $P \in {\PP} 
\setminus \cup_{i \neq j} (C_i \cap C_j)$. 
Thus $(v)_\XX = F$, where  the support $F$ does
not contain any $C_i$ and any point in ${\PP} \setminus \cup_{i
  \neq j} C_i \cap C_j$.   Let $z$ be as in (\ref{saltman1})  and $u =vz$. 
 Since $v  = u_i$ modulo $(\pi_i)$ and 
 $z$ maps to an $l^{\rm th}$ power in $\kappa(C_i)^*$ for each $i$,
$u $ maps to $u_i$ in $\kappa(C_i)^*/\kappa(C_i)^{*l}$ for all $i$. 
In particular $u$ is a unit at $C_i$ for all $i$.  We have 
$(u)_\XX = (z)_\XX + (v)_\XX = E + E'$, where $E' =  Z + F$.
Since the support of $Z$  does not pass through 
any point in $\PP$ and the support of $F$ does not contain any $C_i$ and
does not pass through any point in $\PP \setminus \cup_{i \neq j} (C_i \cap C_j)$, 
$u$ has the  the required properties. 
\end{proof}

\vskip 3mm

The following Lemma is due to Saltman (\cite{S3}, 7.12). We
record here a proof due to Saltman, which bridges missing details in
the original proof.

\begin{lemma}
\label{saltman3} 
 Let   $u_i$, $ 1\leq i \leq n$ and $u$ be as in (\ref{saltman2}). 
 Then there
exists $x \in K^*$ such that $x$ is a 
norm from $K(\sqrt[l]{u})$ and  $(x)_{\XX} = E + E''$ where $E''$
does not pass through  any point  in ${\PP} \cup (\cup_i (E'
\cap C_i))$.   
\end{lemma}

\begin{proof}     Let ${\PP}' = {\PP}  \cup (\cup_i (E'
\cap C_i))$ and  $R$ be the semi-local ring at ${\PP}'$. Let $T$ be the
integral closure of $R$ in $L = K(\sqrt[l]{u})$. 
Since  $L/K$ is ramified at $E_j$ for all $j$, 
there is a unique prime $\tilde{E}_j$  in Spec$(T)$ 
lying over $E_j$.

Since $T$ is normal, it is regular at codimension one points. We
claim that the divisor $\sum t_j \tilde{E}_j$ is principal on $T$.
Since $T$ is semi-local, it is enough to verify that $\sum t_j
\tilde{E}_j$ is principal at each maximal ideal $m$ of $T$. 
 Let
$m$ be a maximal ideal of $T$ and  $P
\in {\PP}'$ be the point corresponding to the maximal ideal of $R$
lying below $m$. If  $P$ is not on any of the $E_j$'s, then clearly
$\sum t_j\tilde{E}_j$ is principal at $m$.

 Suppose that $P$ is on
$E_j$ for some $j$.   Since $E_j$  avoids $C_i \cap C_{i'}$ for
$i \neq i'$, $P \not\in C_i \cap C_{i'}$ for all $i \neq i'$.
Since $E'$ avoids $\PP \setminus \cup_{i \neq i'}(C_i \cap C_{i'})$
and $C_i \cap E_j  \subset \PP \setminus \cup_{i \neq i'}(C_i \cap C_{i'})$, 
$E'$ avoids $C_i \cap E_j$ for all $i$.  Hence $P \not\in E' \cap C_i$ for all $i$. 
 Thus  $P \in \PP'  \setminus  (\cup_i (E' \cap C_i)) =   \PP$.
Since the the support of $E'$ avoids $\PP \setminus \cup_{i\neq i'}(C_i \cap C_{i'})$
and $P \in \PP \setminus \cup_{i \neq i'}(C_i \cap C_{i'})$, 
$P$ is not in the support of $E'$. 
Hence,  the divisor  of $u$ at $P$ is $E  = \sum t_jE_j$. In particular, the divisor of
$\sqrt[l]{u}$ at $m$ is $\sum t_j \tilde{E}_j$.
We conclude that $\sum t_j \tilde{E}_j$ is principal at $m$.
 
Let $w \in T$ be such that the divisor of $(w)_T = \sum t_j
\tilde{E}_j$. Let $x = N_{L/K}(w)$. We claim that $x$ has the
required properties. The divisor of $x$ on $R$ is $\sum t_j E_j$.
Hence $(x)_{\XX}$  is of the form $\sum
t_j E_j + E''$ for some $E''$ with the support of $E''$  avoiding all the points of ${\PP}'$. 
 \end{proof}

\vskip 3mm

Let ${\XX}$, $K$, $\{ C_1, \cdots , C_n \}$, $E$, ${\PP}$ and
$B$ be as above.  Let $\alpha \in H^2(K, \mu_l)$ be  a symbol.  Assume that
there are no chilly loops  for $\alpha$ on ${\XX}$ 
and  $\rama \subset \{ C_1, \cdots, C_n \}$. Since $\alpha$
is a symbol, there are no hot points for $\alpha$ on ${\XX}$ (\cite{S2},
2.5).
By (\cite{S3}, 7.10), for every $C_i$ in 
$\rama$, there exists $u_i \in B/(\pi_i)$  such that $u_i$ is a unit
at all $P \in {\PP} \cap C_i$ except possibly at $C_i \cap C_j$ for
$i \neq j$ and $\partial_{C_i}(\alpha) = (u_i)$.  For each $C_j$ not
in $\rama$, suppose  $u_j \in B/(\pi_j)$ is  a non-zero element  such
that $u_j(P) = u_i(P)$  for all $P \in C_i \cap C_j$ and is a unit at
all $P \in {\PP} \setminus \cup_{i \neq j}(C_i \cap  C_j)$.  Let $u
\in K^*$ be in (\ref{saltman2}) for this  choice of $u_i \in B/(\pi_i)$.

\begin{lemma}
\label{saltman4}
  Let   $u_i$, $ 1\leq i \leq n$, $u$ be as in (\ref{saltman2}) and   $x$ and  $E''$ be as in  (\ref{saltman3}).     
Let $D$ be an irreducible curve in the support of $E''$ with coefficient
prime to $l$. Then  $\alpha$ is unramified at $D$ and  the
specialisation of $\alpha $ at  $D$ is unramified at every discrete
valuation of $\kappa(D)$ centered on a closed point of $D$. 
\end{lemma}

\begin{proof}   Since $E''$ avoids $\PP$ and $\PP$ contains 
at least one point from $C_i$ and $E_j$ for all $i, j$,
$D \neq C_i$ and $D \neq E_j$ for all $i, j$. In particular,
$\alpha$ is unramified at $D$. Let $\overline{\alpha}$ be the specialization of
$\alpha$ at $D$. Let $P$ be a closed point of $D$. If $P $ is not in
$C_i$ for any $C_i$ in $\rama$, then $\overline{\alpha}$ is unramified at every
discrete valuation of $\kappa(D)$ centered on $P$. Suppose that $P \in
C_i \cap D$ for some $C_i$ in $\rama$.   Since $E''$ avoids $\cup_i (E'
\cap C_i) $, $D$ avoids $\cup_i (E'
\cap C_i)$.  If $D$ is in the support of $E'$, then $P \in E'' \cap E' \cap C_i$, 
leading to a contradiction to the fact that $E''$ avoids $E' \cap C_i$. Hence
$D$ is not in the support of $E'$.
 
Let $\YY$ be  the normalization of $\XX$
in $K(\sqrt[l]{u})$. 
 Since supp$_\XX(u)$ is contained in 
 supp$(E) \cup$ supp$(E')$  and $D$  is not in the support of $E$ add $E'$,
   the divisor $D$ is unramified in $\YY$. 
Since $x$ is a norm  from $K(\sqrt[l]{u})$ and the divisor 
$(x)_\XX$ of $x$ on $\XX$ contains $D$ with multiplicity prime to $l$,  
the pull back of $D$ on $\YY$ is of the form $D_1 + \cdot + D_l$ for 
some distinct irreducible curves $D_i$ on $\YY$. 
Hence  the valuation associated to $D$ on $K$ splits in $K(\sqrt[l]{u})$
and the image of $u$ in $\kappa(D)$ is an
$l^{\rm th}$ power. 

We have $(u)_\XX = E + E'$ and 
the support of $E$ and the support of $E'$
do not contain any point in Supp$(E'') \cap C_i$.
Further $P \in $ supp$(E'') \cap C_i$, so that  $u$ is a unit in
$A_P$.   Let  $\delta$ be a prime defining $D$ at $P$. Let $\overline{R}_{\delta}$ be the
integral closure of $A_P/(\delta)$ in $\kappa(D)$. Since $u$ is a unit at $P$ and the image
$\overline{u}$ of $u$ in $\kappa(D)$ is an $l^{\rm th}$ power,   $\overline{u}$ is
an $l^{\rm th}$ power  in $\overline{R}_{\delta}$.

Since  $D$ avoids all the points in ${\PP}$, $P$ is a curve point of $\alpha$. 
Hence $\partial_{C_i}(\alpha)$ is
unramified at $P$,   $\partial_{C_i}(\alpha) = (\overline{u})$ and
 $u$ is a unit at $P$,  We  have $\alpha = \alpha' +  (u) \cdot ( \pi_i)$   
 with  $\alpha'$ unramified at $P$.  Thus the
specialisation  $\overline{\alpha}$ of $\alpha$ at $D$ is equal to
$\overline{\alpha'} + (\overline{u}) \cdot
(\overline{\pi_i})$.  Let $R$ be a discrete valuation ring of
$\kappa(D)$ centered at $P$.   Since $A_P/(\delta)$ is the local ring at
$P$ on $D$, $A_P/(\delta) \subset R$ and
hence   $\overline{R}_{\delta} \subset R$.  Since $\alpha'$ is unramified at $P$, 
$\overline{\alpha'}$ is unramified on $A_P/(\delta)$ and hence $\overline{\alpha'}$
is unramified at $R$. 
Since $\overline{u}$ is an $l^{\rm
  th}$ power   in  $\overline{R}_{\delta}$,  $\overline{u}$ is an $l^{\rm th} $ power in
$R$.  In particular $\overline{\alpha}$ is unramified at $R$.
\end{proof}

\begin{lemma}
\label{saltman5}
   Let $g \in K^*$.  Suppose
$(g)_{\XX}  = \sum r_iC_i + G$ for  some divisor $G$ 
on ${\XX}$ which does not pass through any  point in ${\PP
  }$ and $r_i \geq 0$.   Let $P$ be a closed point on 
$C_i$ for some $i$.  
If   $\partial_{C_i} ( \alpha \cdot (g))$ is unramified at  $P$, 
then either   $\partial_{C_i}(\alpha) \in H^1(\kappa(C_i), \mu_l)$ is
unramified and split at $P$ or the intersection multiplicity $(C_i
\cdot G)_P$ is a multiple of $l$. 
\end{lemma}

\begin{proof}  Let $B$ and $\pi_i$ be as above. 
If $P$ is not on the support of $G$, then $(C_i
\cdot G)_P =0$. Assume that $P $ is in the support of $G$.  Since the
support of $G$ does not contain any point of ${\PP}$,  $P \not\in
C_j$ for any $j \neq i$. Hence $P$ is a curve point of $\alpha$ and
$\alpha =  \alpha'+ (u) \cdot (\pi_i)$ for some $\alpha'$ unramified at
$P$ and $u$ a unit at $P$.   We have $g = \pi_i^{r_i} \theta$ for  some
$\theta \in A_P$ which is not divisible by $\pi_i$. Then  $(C_i \cdot G)_P
= \nu_{P}(\overline{\theta})$,  where $\nu_P$ is the discrete
valuation on $\kappa(C_i)$ at $P$ and  $\overline{\theta}$ is the
image of $\theta$ in $A_P/(\pi_i)$.  
 
We  have $\alpha \cdot (g) = \alpha' \cdot (\pi_i^{r_i} \theta) + (u)
\cdot (\pi_i) \cdot (\pi_i^{r_i}\theta)$ and  $\partial_{C_i}(\alpha
\cdot (g)) =   \overline{\alpha'}^{r_i} + (\overline{u}) \cdot
(\overline{\theta}^{-1}) $, where $\overline{\alpha'}$ is the
specialisation of $\alpha$ at $C_i$.  Since $\alpha'$ is unramified at $P$, we have
$\partial_P(\partial_{C_i}(\alpha \cdot (g))) =
(u(P))^{-\nu_P(\overline{\theta})}$.  If $\partial_{C_i}(\alpha \cdot
(g))$ is unramified at $P$, then $(u(P))^{-\nu_P(\overline{\theta})} =
1 \in H^1(\kappa(P), \mu_l)$. Thus either $u(P)$ is an $l^{\rm th}$
power in $\kappa(P)$ or $\nu_P(\overline{\theta})$ is a multiple of
$l$.  \end{proof}

\begin{remark} 
\label{remark} Let $\alpha$ and $g$ as in (\ref{saltman5}).  Suppose
that  $r_i$ is coprime to $l$. Let $r_i'$ be the inverse of $r_i$
modulo $l$.  Then $\partial_{C_i}(\alpha \cdot (g)) =
r_i'\beta_{C_i}$, where $\beta_{C_i}$
is the residual Brauer class of $\alpha $ with respect to $g$ as defined in (\cite{S2}, p.837). 
 If $(r_i, l ) = 1$,  
the above lemma is the same as (\cite{S2}, 4.2(c)). In particular, we have a new interpretation of
the residual Brauer class. 
\end{remark}

\section{A local-global principle}

Let $\XX$ be an excellent  regular integral surface which is
quasi-projective over an affine scheme. Let   $K$ be the   function
field of ${\XX}$. Let $l$ be a prime which is a unit on ${\XX} $.  
We say that $K$ satisfies  {\it local-global principle for $H^3(K,
~\mu_l^{\otimes 3})$ in terms of symbols in $H^2(K, ~\mu_l^{\otimes 2})$} if the following
holds:  for any  $\zeta \in H^3(K, ~\mu_l^{\otimes 3})$ and a symbol $\alpha \in
H^2(K, ~\mu_l^{\otimes 2})$, if for any discrete valuation  $\nu  $ of $K$ there exists $f_{\nu} \in
K_{\nu}^*$ such that $\zeta - \alpha \cdot (f_{\nu}) \in H^3_{nr}(K_{\nu},
~\mu_l^{\otimes 3})$, then there exists $ f\in K^*$ such that $\zeta - \alpha
\cdot (f) \in H^3_{nr}(K, ~\mu_l^{\otimes 3})$.

In this section we prove that $K$ satisfies local-global principle for
$H^3(K, \mu_l^{\otimes 3})$ in terms of symbols in $H^2(K, \mu_l^{\otimes 2})$ under some
conditions on $K$.

 The following notation and assumptions  are  used in (\ref{choice-at-chilly1},
\ref{choice-at-chilly2}, \ref{choice-local}).\\
$\bullet$ $\XX$ an excellent regular surface, quasi-projective over an affine scheme\\
 $\bullet$ $K$ contains a primitive $l^{\rm th}$ root of unity. \\
 $\bullet$ $\alpha \in H^2(K, \mu_l)$ a symbol \\
$\bullet$  $\zeta \in H^3(K, \mu_l)$ \\
$\bullet$ $\rama \cup \ramz$  is a 
union of regular curves  with only normal crossings \\
$\bullet$ $\rama = \{ C_1, \cdots , C_m \}$ and $\rama \cup \ramz = $  $\{ C_1,
\cdots , C_m,   C_{m+1}, \cdots , C_n \}$.

For each closed point $P$ of ${\XX}$, let $A_P$ be the local ring at $P$
and $\hat{A}_P$ be the completion of $A_P$ at its maximal ideal $m_p$.

Let ${\PP}$ be a finite set of closed points containing all the
nodal points of $\rama \cup \ramz$ and  at least one point from each
$C_i$.    Let  $B$ be the semi-local ring at  ${\PP
 }$. Let $\pi_i \in B$ be a prime defining $C_i$ for $1 \leq i \leq
n$.

\begin{lemma}
\label{choice-at-chilly1}
 Suppose  that  $H^2(\kappa(P), \mu_l) = 0$  for every closed point of  $P$ of ${\XX}$
 and  for each $x \in {\XX}^1$ there exists $f_x \in K_x^*$ such that 
$\zeta - \alpha \cdot (f_x) \in H^3_{nr}(K_x, \mu_l)$.   
 Let $P \in C_i \cap C_j$ be a chilly point of
$\alpha$.  Given any integer  $r_i$  with  $0 \leq r_i \leq l-1$,  there
exists  an integer  $r_j$ such that $0 \leq r_j \leq   l-1$
and  $\zeta - \alpha \cdot (\pi_i^{r_i}\pi_j^{r_j})$ is unramified on $\hat{A}_P$. 
\end{lemma}  

\begin{proof}  By the assumption, there exists $ f_i \in
K_{C_i}^*$ such that $\zeta - \alpha \cdot (f_i) \in H^3_{nr}(K_{C_i},
\mu_l)$.  In particular, $\partial_{C_i}(\zeta)
= \partial_{C_i}(\alpha \cdot (f_i)) \in H^2(\kappa(C_i), \mu_l)$.
Since $P$ is a chilly point of $\alpha$, we have $\alpha = (u) \cdot
(\pi_i\pi_j^{s_j}) \in H^2(K_P, \mu_l)$ (cf. \ref{saltman-2dim}) for some unit $u$ in
$A_P$ and $1 \leq s_j \leq l-1$.  Then 
$$
\begin{array}{lcl}
\partial_{C_i}(\zeta) & = & \partial_{C_i}(\alpha \cdot (f_i)) \\
& = &  \partial_{C_i}((u) \cdot (\pi_i\pi_j^{s_j}) \cdot (f_i)) \\ 
& =  &  (\overline{u}) \cdot \partial_{C_i}((\pi_i\pi_j^{s_j}) \cdot (f_i)) \\
& = & (\overline{u}) \cdot  (g_i) \in H^2(\kappa(C_i)_P, \mu_l)
\end{array}
$$
 for some
$g_i \in \kappa(C_i)^*$, where $\overline{u}$ is the image of $u$ in
$\kappa(C_i)$.  Since $\kappa(C_i)$ is the field of fractions of
$A_P/(\pi_i)$ and $A_P/(\pi_i)$ is a discrete valuation ring with
$\overline{\pi}_j$ as a parameter, we have $g_i =
w_i\overline{\pi}_j^{t_i}$ for some unit $w_i \in A_P/(\pi_i)$ and
$t_i \geq 0$.  Hence $\partial_P(\partial_{C_i}(\zeta))
= \partial_P((\overline{u}) \cdot (g_i)) = (u(P))^{t_i}$.  The point
$P$ being a chilly point of $\alpha$, $u(P)$ is not an $l^{th}$ power
in $\kappa(P)$ and $0 \leq t_i \leq l-1$.  Let $r_j $
be such that $0 \leq r_j \leq l-1$ and $r_j = r_is_j - t_i$ modulo
$l$.
 
We now show that $\zeta - \alpha \cdot (\pi_i^{r_i}\pi_j^{r_j})$ is
unramified at $\hat{A}_P$. Since $\zeta$ and $\alpha \cdot
(\pi_i^{r_i}\pi_j^{r_j})$ are ramified on $\hat{A}_P$ at most at
$\pi_i$ and $\pi_j$, it is enough to show that the residues of $\zeta-
\alpha \cdot (\pi_i^{r_i}\pi_j^{r_j})$ are equal to zero at $\pi_i$ and
$\pi_j$.    

 We have $\partial_{C_i}(\alpha \cdot (\pi_i^{r_i}\pi_j^{r_j}))
= \partial_{C_i}((u) \cdot (\pi_i \pi_j^{s_j}) \cdot
(\pi_i^{r_i}\pi_j^{r_j})) = (\overline{u}) \cdot (
\overline{\pi}_j^{r_is_j - r_j}) = (\overline{u}) \cdot (
\overline{\pi}_j^{t_i}) \in H^2(\kappa(C_i)_P, \mu_l)$.  Since
$\hat{A}_P/(\pi_i)$ is a complete discrete valuation ring with residue
field $\kappa(P)$ and $H^2(\kappa(P), \mu_l) = 0$, $(\overline{u})
\cdot (w_i) = 0 \in H^2(\kappa(C_i)_P, \mu_l)$ and $(\overline{u})
\cdot (w_i\overline{\pi}_j^{t_i}) = (\overline{u}) \cdot
(\overline{\pi}_j^{t_i}) \in H^2(\kappa(C_i)_P, \mu_l)$.  Thus
$\partial_{C_i}(\alpha \cdot (\pi_i^{r_i}\pi_j^{r_j}))  = 
(\overline{u}) \cdot (w_i\overline{\pi}_j^{t_j}) =
 (\overline{u}) \cdot (g_i) =
\partial_{C_i}(\zeta)$
in $H^2(\kappa(C_i)_P, \mu_l)$.  Hence $\zeta - \alpha \cdot
(\pi_i^{r_i}\pi_j^{r_j})$ is unramified on $\hat{A}_P$ except possibly
at $(\pi_j)$. Thus,  by (\ref{ramified-at-pi}), $\zeta-\alpha \cdot (\pi_i^{r_i}\pi_j^{r_j})$
is unramified at $\hat{A}_P$. 
\end{proof}

\begin{lemma}  
\label{choice-at-chilly2} 
Suppose  that $H^2(\kappa(P), \mu_l)  = 0$ and  there are no chilly loops
for $\alpha$ on ${\XX}$.  Then there exist integers $r_i$, $ i =
1, \cdots   n$ such that  for each   point $P \in C_i \cap
C_j$ $(i \neq j)$ there exists a unit $w_P \in A_P$ such that    
$\zeta  - \alpha \cdot  (w_P \pi_i^{r_i} \pi_j^{r_j})$  is
unramified on $\hat{A}_P$.    
\end{lemma}

\begin{proof} For each $i$, $ 1 \leq i \leq n$, 
there exists $f_i \in K_{C_i}^*$ such that $\zeta - (\alpha) \cdot (f_i)$ is
unramified at $C_i$.  

Suppose $m+1 \leq i \leq n$. 
Since $\alpha$ is unramified at $C_i$ and
$\zeta$ is ramified at $C_i$,  we have $\partial_{C_i}(\zeta) =
\overline{\alpha}^{\nu_{C_i}(f_i)}$, where $\overline{\alpha} \in
H^2(\kappa(C_i), \mu_l)$ is the specialisation of $\alpha$ at $C_i$
and $\nu_{C_i}$ is the valuation at $C_i$.
Since $\zeta$ is ramified at $C_i$, $\overline{\alpha}  \neq 0$ and   
  $\nu_{C_i}(f_i)$  modulo $l$  is independent of the choice of
$f_i$.  Let $r_i = \nu_{C_i}(f_i)$. 

Let $1 \leq i \leq m$.   If $C_i$ has no chilly point, then we choose $r_i
= 1$.  Since there are no chilly loops, using (\ref{choice-at-chilly1}), we construct $r_i$
corresponding to each $C_i$ which has a chilly point on it,   
such that $\zeta -  \alpha \cdot
(\pi_i^{r_i}\pi_j^{r_j})$ is  unramified on  $\hat{A}_P$ for all
chilly points $P$ of $\alpha$.  

We now show that these $r_i$'s have the required property for suitable choices of 
$w_P $  as $P$ various over points in $C_i \cap C_j$ for $ 1 \leq i < j \leq n$. 
Let $P \in C_i \cap C_j$ for some $1 \leq i < j \leq n$.  

Suppose that $j \leq m$.
If $P$ is a chilly point, then $\zeta - \alpha \cdot
(\pi_i^{r_i}\pi_j^{r_j})$ is unramified at $\hat{A}_P$ by the choice of $r_i$. In
this case we set $w_P = 1$.
 
Suppose that $P$ is a cold point.   Then
$\alpha =    (u\pi_j) \cdot (v \pi_i^s) \in H^2(K_P, \mu_l)$
 for some  units $u$, $v$ at $P$ and $1
\leq s \leq l-1$ (cf., \ref{saltman-2dim}).    Since
$\hat{A}_P/(\pi_i)$ is a  complete discrete valuation ring with
$\overline{\pi_j}$ as a  parameter and   $H^2(\kappa(P), \mu_l) = 0$,
  every element of $H^2(\kappa(C_i)_P, \mu_l)$ is of the form 
$ (\overline{\pi}_j) \cdot (w) $ for some unit $w$ in $\hat{A}_P/(\pi_i)$.
Hence we have  $\partial_{C_i}(\zeta - \alpha \cdot (\pi_i^{r_i}
\pi_j^{r_j}))  =  (\overline{\pi_j}) \cdot (w_i) 
\in H^2(\kappa(C_i)_P, \mu_l)$ for some unit $w_i \in \hat{A}_P /(\pi_i)$.  
Since $s$ is coprime to $l$, there exists a unit  $w_P \in A_P$  
such that $w_P^{-s}(P) = w_i(P)$ modulo $\kappa(P)^{*l}$.
In particular $w_P^{-s} = w_i$ modulo $K_P^{*l}$. 
Over $K_P$, we have  
$$
\begin{array}{rcl}
\partial_{C_i}(\zeta -  \alpha \cdot
(w_P\pi_i^{r_i}\pi_j^{r_j})) & = & \partial_{C_i}(\zeta -  \alpha \cdot
(\pi_i^{r_i}\pi_j^{r_j}) ) - \partial_{C_i}(\alpha \cdot
(w_P) ) \\
& = & (\overline{\pi_j}) \cdot (w_i)  -  \partial_{C_i}((u\pi_j) \cdot (v\pi_i^s) \cdot (w_P)) \\
& = & (\overline{\pi_j}) \cdot (w_i)  -  (\overline{u}\overline{\pi}_j) \cdot 
 \partial_{C_i}( (v\pi_i^s) \cdot (w_P)) \\
 & = & (\overline{\pi_j}) \cdot (w_i)  -  (\overline{u}\overline{\pi}_j) \cdot (w_P^{-s}) \\ 
& = & (\overline{\pi_j}) \cdot (w_i)  -  (\overline{u}\overline{\pi}_j) \cdot (w_i) \\
& = & (\overline{\pi_j}) \cdot (w_i)  -  (\overline{\pi}_j) \cdot (w_i)  - (\overline{u} ) \cdot (w_i) \\
& = & -(\overline{u}) \cdot (w_i)  \\
& = & 0.
\end{array}
 $$
  In particular, $\zeta
- \alpha \cdot (w_P \pi_i^{r_i}\pi_j^{r_j})$ is unramified on
$\hat{A}_P$ except possibly at $C_j$. By (\ref{ramified-at-pi}), $\zeta - \alpha \cdot
(w_P\pi_i^{r_i}\pi_j^{r_j})$ is unramified on $\hat{A}_P$.  

Suppose
that $P$ is a cool point of $\alpha$. Then $\alpha =  0 \in H^2(K_P,
\mu_l)$ (cf. \ref{saltman-2dim}) and hence $\partial_P(\partial_{C_i}(
\alpha \cdot (f_i)) = 0$ (cf. the diagram in the proof of
\ref{residue-unramified-at-P}).  Since $\zeta - \alpha \cdot (f_i)$ is unramified at $C_i$,
 we have $\partial_P(\partial_{C_i}(\zeta))  =
\partial_P(\partial_{C_i}(\alpha \cdot (f_i))) = 0$.  Since $\kappa(C_i)_P $
is a complete discrete valued field with residue field $\kappa(P)$ and
$H^2(\kappa(P), \mu_l) = 0$,  $\partial_{C_i}(\zeta) $ is zero in $H^2(\kappa(C_i)_P,
\mu_l)$.  By (\ref{ramified-at-pi}), $\zeta$ is unramified  on $\hat{A}_P$ and
hence $\zeta -  \alpha \cdot (f_P)$ is unramified on
$\hat{A}_P$ for any $f_P \in K_P$.  We set $w_P = 1$. Then $\zeta - \alpha \cdot
(w_P\pi_i^{r_i}\pi_j^{r_j})$ is unramified on $\hat{A}_P$.

Suppose $j \geq m+1$.  Then $\alpha$ is unramified at $C_j$. Hence, by
the choice of $r_j$,  $\partial_{C_j}(\zeta) =
\overline{\alpha}^{r_j}$, where $\overline{\alpha}$ is the
specialisation of $\alpha$ at $C_j$.  We have  $\partial_{C_j}(\alpha
\cdot (\pi_i^{r_i}\pi_j^{r_j})) = \overline{\alpha}^{r_j}$. Hence
$\partial_{C_j}(\zeta - \alpha \cdot (\pi_i^{r_i} \pi_j^{r_j})) = 0$.
Thus, by (\ref{ramified-at-pi}),  $ \zeta - \alpha \cdot (w_p\pi_i^{r_i} \pi_j^{r_j})$ is
unramified on $\hat{A}_P$ with $w_P = 1$.   \end{proof}

\begin{lemma}
\label{choice-local} Suppose  that there are no chilly loops
for $\alpha$ on ${\XX}$. 
    Let $r_i$, $1 \leq i \leq n$ be 
as in (\ref{choice-at-chilly2}).
 There exists $f \in K^*$ such that $(f)_{\XX
  } = \sum_{i = 1}^n r_i C_i   + F$  for some divisor $F$ which does
not pass through any point of ${\PP}$   and  
$\zeta - \alpha \cdot (f) $ is unramified at $\hat{A}_P$ for all $P \in {\PP}$. 
\end{lemma}

\begin{proof}  Let $B$ and
$\pi_i$ be as above.  Let $r_i$ be as in (\ref{choice-at-chilly2}).   Let $g  = \prod_1^n
\pi_i^{r_i} \in B$. 

 Let $P \in {\PP}$.   Suppose  $P  \in C_i \cap C_j$ for some $i
 \neq j$.  Since $\{C_1, \cdots , C_n \}$ have   only normal crossings,
we have $g = u_P \pi_i^{r_i} \pi_j^{r_j}$ for some  unit $u_P$ in $A_P$.
Suppose   $P \in  C_i$ for only one $i$.  Then we have $g =
u_P\pi_i^{r_i}$ for some unit $u_P$ in $A_P$.  If $P \not\in C_i$ for
all $i$, let $u_P = 1$.

Let $P \in {\PP}$. If $P \in C_i \cap C_j$ for some $1 \leq i < j
\leq n$, let $w_P$ be as in (\ref{choice-at-chilly2}). If $P \not\in C_i \cap C_j$ for all
$i \neq j$, let
$w_P = 1$. Let $w  \in B$  be a unit such that $w(P) = w_P(P)/u_P(P)$  for
all  $P \in  {\PP}$.  Then $wu_Pw_P^{-1}$ is an $l^{\rm th}$ power
in $K_P$ for all $P \in \PP$. Let $f = wg$.  Then $(f)_{\XX} = \sum_{i=1}^n r_iC_i +  F$  
for some divisor $F$ which does not pass through any point of ${\PP}$.

Let $P \in {\PP}$. If  $P  \in  C_i \cap C_j$ for
some $i \neq j$,   we have $f =  wg = wu_P \pi_j^{r_i}\pi_j^{r_j} = 
w_P\pi_i^{r_i}\pi_j^{r_j}$ modulo
$l^{\rm th}$  powers in $K_P^*$ and $\zeta - \alpha \cdot (f)$ is unramified on $\hat{A}_P$ 
by (\ref{choice-at-chilly2}). Suppose    $P \in C_i$ for some $i$ and
not in $C_i \cap C_j$  for any $i \neq j$.   Then $f = w_P
\pi_i^{r_i}$ and  $\zeta$,  $\alpha$, $\alpha \cdot (f)$ 
 are ramified at most along $C_i$. By (\ref{ramified-at-pi}),
$\zeta$ and $\alpha \cdot (f)$ are  unramified on $\hat{A}_P$.
If  $P \not\in C_i$ for all $i$, then $f$ is a unit at
$P$ and $ \zeta$ and $\alpha \cdot (f)$ are unramified on $\hat{A}_P$. Hence $f$ has all 
the required properties.  \end{proof}

\vskip 3mm

Let $ D$ be an integral scheme of dimension one and $\kappa(D)$ its
function field.  Let $l$ be a prime which is a unit on $D$.  Let $L$
be a finite extension of $\kappa(D)$.  An element $\alpha \in H^n(L,
\mu_l)$ is said to be {\it unramified} over $D$ if $\alpha$ is
unramified at every discrete valuation of $L$ with restriction to
$\kappa(D)$ centered on a closed point of $D$.  Let $H^n_{nr}(L/D,
\mu_l)$ denote the set of elements of $H^n(L, \mu_l)$ which are
unramified over $D$.  

Let ${\XX}$ be an excellent regular integral surface which is
quasi-projective over an affine scheme. Let $l$ be a prime which is a
unit on ${\XX}$. We say that ${\XX}$ is {\it l-special} if it
satisfies the following: \\
\noindent i) $H^2_{nr}(L/D, \mu_l) = 0$ for every irreducible curve
$D$ on ${\XX}$ and for every finite extension $L$ of
$\kappa(D)$. \\ ii) the $l$-cohomological dimension of $\kappa(P)$ is
at most 1 for every closed point $P$ of ${\XX}$.

Let $R$ be the ring of integers in a $p$-adic field and $l $ a prime
not equal to $p$.  Let ${\XX}$ be a two-dimensional excellent
regular integral proper scheme over Spec$(R)$.  Then ${\XX}$ is
$l$-special. A two-dimensional regular integral proper variety over a
finite field of characteristic not equal to $l$ is also
$l$-special. Note that if ${\XX}$ is $l$-special, then any blow up
of ${\XX}$ at a closed point is also $l$-special.

Let ${\XX}$ be an $l$-special scheme.  Let $P$ be a closed point of
$\XX$.  Since the $l$-cohomological dimension of $\kappa(P)$ is at
most 1, $H^n(\kappa(P), \mu_l) = 0$ for $n \geq 2$.  In particular, if
$K_P$ is the field of fractions of the completion of ${\OO}_{{\XX
}, P}$ at the maximal ideal, $H^n_{nr}(K_P, \mu_l) = 0$ for $n \geq
2$.  We use this fact throughout the   rest of the paper.

\begin{theorem}
\label{main}
 Let $l$ be a prime. Let ${\XX}$ be an
$l$-special scheme and $K$ its function field.  Suppose that $K$
contains a primitive $l^{\rm th}$ of unity.  Let $\alpha \in H^2(K,
\mu_l)$ be a symbol and $\zeta \in H^3(K, \mu_l)$.  Suppose that
$\rama \cup \ramz$ is a union of regular curves with normal crossings
on ${\XX}$.  Suppose for every $x \in {\XX}^1$, there exists
$f_x \in K^*_x$ such that $\zeta - \alpha \cdot (f_x) \in
H^3_{nr}(K_x, \mu_l)$.  Then there exists $f \in K^*$ such that $\zeta
- \alpha \cdot (f) \in H^3_{nr}(K/{\XX}, \mu_l)$.
\end{theorem}

\begin{proof} Let $\rama = \{ C_1, \cdots, C_m \}$ and  $\rama \cup \ramz = \{ C_1,
\cdots , C_m, C_{m+1}, \\ \cdots  , C_n \}$. Let ${\PP}$ be a finite set of closed 
points containing $C_i\cap
C_j$ for $1 \leq i < j \leq n$ and at least one point from each
$C_i$. Let $B$ be the semi-local regular ring at all points in ${\PP
 }$. Let $\pi_i \in B$ be a prime defining $C_i$ on $B$.

\vskip 3mm

\noindent
{\it Step 1:   Reduction to the case where 
there are no chilly loops on $\XX$ for $\alpha$:}

  Let $P \in C_i \cap C_j$ be
a chilly point of $\alpha$.  Then, by (\ref{choice-at-chilly1}), taking $r_i 
= 0$, there exists $r_j$
such that $\zeta - \alpha \cdot ( \pi_j^{r_j})$ is
unramified on $\hat{A}_P$. Let ${\XX}'$ be the  simple  blow-up of ${\XX}$
at $P$ and $\tilde{E}$ be the  exceptional curve on ${\XX}'$.   
 Since $\zeta - \alpha \cdot
(\pi_j^{r_j})$ is unramified on $\hat{A}_P$ and $H^2(\kappa(P), \mu_l) =0$,
by (\ref{purity}),  $\zeta - \alpha \cdot
(\pi_j^{r_j})$  is 
unramified at $\tilde{E}$.   Hence, by blowing-up chilly points we
still have the hypothesis of the theorem. We therefore assume that
there no chilly loops on ${\XX}$ for $\alpha$ (\cite{S2},  Corollary 2.9). 

\vskip 3mm

\noindent
{\it Step 2:  We show that there  exists $g \in K^*$  satisfying \\
\noindent
i)  $\zeta  - \alpha \cdot (g)$ is unramified at all $C_i$ \\
\noindent
ii)  $(g)_\XX = \sum r_i C_i + E$  with  $ r_i \in \Z$,  $E$  a divisor which does
not pass through any point in $\PP$ \\
\noindent
iii)  $(C_i \cdot E)_P$ is a  multiple of $l$ for all $P \in C_i \cap E$ for all $i$.}

Since there are no chilly loops on $\XX$ for $\alpha$, 
choose $r_i$ as in (\ref{choice-at-chilly2}). 
 Using (\ref{choice-local}), we choose  $ h \in K^*$ such that  
 $\zeta - \alpha \cdot (h)$ is unramified on $\hat{A}_P$ for all $P \in {\PP}$
 and $(h)_\XX = \sum_i r_iC_i + F$ with $F$ a divisor on $\XX$ which 
 does not pass through any point in $\PP$.

 Since $\alpha$ is a symbol, we have  $\alpha
= (u_i) \cdot (v_i\pi_i^{\epsilon})$  for some units $u_i$ and $v_i$
at $C_i$ and $\epsilon = 0, 1$ (cf. \ref{dvr-symbol}).   
For any $a \in K^*$, we have 
$$
\partial_{C_i}(\alpha \cdot (a)) = \partial_{C_i}((u_i) \cdot (v_i \pi^{\epsilon}) \cdot (a))
 = (\overline{u_i}) \cdot \partial_{C_i}((v_i \pi^{\epsilon}) \cdot (a)) = 
 (\overline{u_i})  \cdot (c_a)
 $$
 for some $c_a \in \kappa(C_i)^*$.
 In particular $\partial_{C_i}(\alpha \cdot (h) ) = (\overline{u}_i) \cdot (c_i)$ for
 some $c_i \in \kappa(C_i)$. 
 Since 
$\zeta - \alpha  \cdot (f_i)$ is unramified at $C_i$ for some $f_i \in
K^*$, we have 
$$\partial_{C_i}(\zeta) = \partial_{C_i}(\alpha \cdot
(f_i)) = (\overline{u_i}) \cdot (a_i)$$ 
for some $a_i \in
\kappa(C_i)^*$.     
Set $\beta_i = \partial_{C_i}(\zeta - \alpha \cdot (h)) \in
H^2(\kappa(C_i), \mu_l)$.
Then 
$$\beta_i = \partial_{C_i}(\zeta ) - \partial_{C_i}(\alpha \cdot (h)) =
 (\overline{u}_i) \cdot (a_i) - (\overline{u}_i) \cdot (c_i) =   
 (\overline{u}_i) \cdot (b_i)$$
 with  $b_i = a_ic_i \in \kappa(C_i)^*$.

Let $P \in {\PP }
\cap C_i$.   By the choice of $h$,   $\zeta - \alpha \cdot (h)$ is unramified
at $\hat{A}_P$. Thus $\beta_i = \partial_{C_i}(\zeta - \alpha \cdot (h))$ is
unramified at $P$ (cf. \ref{residue-unramified-at-P}).  
 Since $H^2(\kappa(P), \mu_l) = 0$ and   $\beta_i$ unramified at
$P$, $\beta_i = (\overline{u}_i) \cdot (b_i)$ is zero over the completion $\kappa(C_i)_P$  of
$\kappa(C_i)$ at  $P$.
Hence $b_i$ is a norm   from the extension  $\kappa(C_i)( \sqrt[l]{\overline{u}_i}) \otimes \kappa(C_i)_P$ 
over  $\kappa(C_i)_P$. Let $M_i = \kappa(C_i)(\sqrt[l]{\overline{u}_i})$.
Let $\theta_P \in M_i  \otimes \kappa(C_i)_P $ be
such that 
$$N_{M_i \otimes \kappa(C_i)_P/\kappa(C_i)_P)}(\theta_P) = b_i.$$

Note that  $B/(\pi_i)$ is a semi-local ring with field of fractions
$\kappa(C_i)$ with maximal ideals corresponding to the points of
${\PP} \cap C_i$.   By the weak approximation,  we choose
$\theta \in M_i$  such that 
$\theta$ is ``sufficiently close'' to $\theta_P$ for all $ P \in C_i \cap \PP$.
Let 
$$ a= N_{M_i/\kappa(C_i)}(\theta).$$
Then $(ab_i)(P) = 1$ for all $P \in C_i \cap \PP$ and  $(\overline{u}_i) \cdot  (a) = 0$.
In particualr,  $(\overline{u}_i) \cdot 
(b_i)    = (\overline{u}_i) \cdot (ab_i)$. Thus,  replacing  $b_i$ by $ab_i$,
 we  assume that $b_i \in B/(\pi_i)$ and $b_i(P) = 1$ for all $P \in C_i \cap {\PP}$.
 By, (\cite{S3}, Proposition 0.3),  there exists a unit $b \in B$   such that $b = b_i \in B/(\pi_i)$ 
 for all $i$.   We set $g = bh$ and prove that $g$ has the required properties $i)$, $ii)$ and $iii)$. 
 Since $b$ is a unit in $B$, the divisor $(b)_\XX$ does not pass through any point in $\PP$. 
By the choice of $h$, we have $(g)_{\XX
  } = (bh)_\XX = \sum r_i C_i +  E$ for some divisor $E$ on ${\XX}$ which
does not pass through any point in ${\PP}$.  Since $b$ is 
a unit at $C_i$,  $\partial_{C_i}(\alpha \cdot (b)) = \partial_{C_i}(\alpha) \cdot (b_i) =
(\overline{u}_i) \cdot (b_i)$. Hence 
$\partial_{C_i}(\zeta - \alpha \cdot (g)) = \partial_{C_i}(\zeta
-\alpha \cdot (h)) - \partial_{C_i}(\alpha \cdot (b)) = \beta_i -
(\overline{u}_i) \cdot (b_i) = 0$ for all $i$.

Let $P \in C_i \cap E$. Then $P \not\in C_i \cap C_j$ for all $i \neq
j$.  Then, by (\ref{ramified-at-pi}), $\zeta$ is unramified on $\hat{A}_P$ and hence
$\partial_{C_i}(\zeta)$ is unramified at $P$ (cf. \ref{residue-unramified-at-P}). Thus
$\partial_{C_i}(\alpha \cdot (g)) = \partial_{C_i}(\zeta)$ 
is unramified at $P$. By (\ref{saltman5}),
either $ \partial_{C_i}(\alpha)$ is unramified and split at $P$ or
$(C_i \cdot E)_P$ is a multiple of $l$.  Suppose that
$\partial_{C_i}(\alpha)$ is unramified and split at $P$. Then $\alpha$
is unramified on $\hat{A}_P$ at $C_i$. Since $P \not \in C_j$ for all
$j \neq i$, $\alpha$ is unramified on $\hat{A}_P$.  Let ${\XX}'$ be
a simple blow-up at $P$ and $\tilde{E}$ the exceptional curve on
${\XX}'$.  Since $\alpha$ is unramified on $\hat{A}_P$, 
$\alpha$ is unramified at $\tilde{E}$ (cf. \ref{purity}).  
Hence $\tilde{E}$ is not in the 
ramification locus of $\alpha$. 
Since $P \not\in C_i \cap C_j$ for all $i \neq j$, 
by (\ref{ramified-at-pi}),  $\zeta$ is unramified on $\hat{A}_P$
and hence $\tilde{E}$ is not in the
ramification locus of $\zeta$ (\ref{purity}). Since ${\XX}'$ satisfies hypothesis
of the theorem, by replacing ${\XX}$ with a blow-up, we assume that $(C_i
\cdot E)_P$ is a multiple of $l$ for all $P \in C_i \cap E$ for all $i$.

\vskip 3mm

\noindent
{\it Step 3: The final choice.}

For $1 \leq i \leq m$, $\partial_{C_i}(\alpha) = (u_i)$. By (\cite{S3},
7.10), we can assume that $u_i(P) = u_j(P)$ for all $P \in C_i \cap
C_j$, $ 1 \leq i < j \leq m$.  Let $P \in {\PP}$. If $P \not\in
C_i$ for $1 \leq i \leq m$, let $u_P$ be an element of $\kappa(P)$ which is
not an $l^{\rm th}$ power.
 If $P \in C_i$ for some $1 \leq i \leq m$, let $u_P =
u_i(P)$, noting that $u_i(P) = u_j(P)$ whenever $P \in C_i \cap C_j$,
$1 \leq i<j \leq m$. For $m+1 \leq j \leq n$,  by Chinese 
remainder theorem, there exist  $u_j \in
\kappa(C_j)^*$  such that $u_j(P) = u_P$ for all $P \in {\PP}
\cap C_j$. For these choice of $u_i$, we choose  $u$  as in
(\ref{saltman2});  namely $(u)_\XX = E + E'$
where $E'$ is a divisor on $\XX$ with the support of $E'$ not passing through 
any point in $\PP \setminus \cup_{i \neq j}{C_i \cap C_j}$ , no $C_i$ or
$E_j$ is in the support of $E'$ and the image of $u$ is equal to $u_i$ in 
$\kappa(C_i)^*/\kappa(C_i)^{*l}$ for all $i$. Let $x$ be a 
norm from $K(\sqrt[l]{u})$ as in (\ref{saltman3}); namely 
$(x)_{\XX} = E + E''$ where $E''$
does not pass through  any point  in ${\PP} \cup (\cup_i (E'
\cap C_i))$.  Let $g$ be as in Step 2) and $f = x^{-1}g$. 
We now show that $\zeta - \alpha \cdot (f)$ is unramified at every
codimension one point of ${\XX}$.  

Let $D$ be a codimension one point of  ${\XX}$.

Suppose $D$ is not in $ram_{\XX}(\alpha) \cup ram_{\XX}(\zeta)
\cup Supp(f)$. Then $\zeta$ and $\alpha \cdot
(f)$ are  unramified at $D$.

Suppose that $D$ is in $ram_{\XX}(\alpha) \cup ram_{\XX
}(\zeta)$. Then $D = C_i$ for some $i$.  By the choice of $g$,
$\zeta - \alpha \cdot (g)$ is unramified at $C_i$ (step 2). 
 Thus it is enough
to show that $\alpha \cdot (x^{-1})$ is unramified at $C_i$.  Since
the image of $u$ in $\kappa(C_i)$ is $u_i$ and $\partial_{C_i}(\alpha)
= (u_i)$, we have $\alpha - (u) \cdot (\pi_i)$ is unramified at $C_i$.
Hence $\alpha \cdot (x^{-1}) - (u) \cdot (\pi_i) \cdot(x^{-1})$ is
unramified at $C_i$.  Suppose that $1 \leq i \leq m$. Since $x$ is a
norm from $K(\sqrt[l]{u})$, $(u) \cdot (\pi_i) \cdot (x^{-1}) = 0$ and
hence $\alpha \cdot (x^{-1})$ is unramified at $C_i$.  Suppose that $i
\geq m+1$. Since $\alpha$ is unramified at $C_i$ and $x$ is a unit at
$C_i$, it follows that $\alpha \cdot (x^{-1})$ is unramified at $C_i$.
 
Suppose that $D$ is in the support of $(f)_{\XX}$
and not in $ram_{\XX}(\alpha) \cup ram_{\XX}(\zeta)$. Then
$\alpha$ and $\zeta$ are unramified at $D$.  If  $\nu_D(f)$ is a
multiple of $l$, then $\zeta - \alpha \cdot (f)$ is unramified at
$D$. Suppose that $\nu_D(f)$ is not a multiple of $l$. 
Since $(g)_\XX = \sum r_iC_i + E$ and $(x) = E + E''$,  we have
$(f)_\XX = (x^{-1}g)_{\XX} = \sum r_iC_i - E'' $. 
Since $D \neq C_i$ for all $i$, $D$ is in the
support   of $E''$. Since $E$ and $E''$ have no common component,
 $\nu_D(g) = 0$  and $\nu_D(f) = \nu_D(x^{-1}g) = 
\nu_D(x^{-1})$  is coprime to $l$.  By (\ref{saltman4}),  
the specialization $\overline{\alpha}$
of $\alpha$ at $D$ is unramified at every discrete valuation of
$\kappa(D)$ centered at any $P \in {\PP} \cap D$. Since $\alpha$ is
ramified only along $C_i$, $1 \leq   i \leq m$, $\overline{\alpha}$ is
unramified at all discrete valuations of $\kappa(D)$ centered at $P$
for all $P \not \in {\PP} \cap D$. Hence $\overline{\alpha}$ is
unramified at every discrete valuation of $\kappa(D)$ centered  on a
closed point of $D$. Since
$H^2_{nr}(\kappa(D), \mu_l) = 0$,  $\overline{\alpha} = 0$.  
Since $\partial_D(\alpha \cdot (f)) =
\overline{\alpha}^{\nu_D(f)} = 0$, $\alpha \cdot (f)$ is
unramified at $D$.  Since $\zeta$ is unramified at $D$, $\zeta -
\alpha \cdot (f)$ is unramified at $D$.  \end{proof}

\begin{cor} 
\label{local-global-ps1}
(cf. \cite{PS2}, 3.4) Let $k$ be a $p$-adic field
and $l$ a prime not equal to $p$. Let $X$ be a curve over $k$ and
$K$ its function field. Assume that $K$ contains a primitive
$l^{\rm th}$ root of unity. Let ${\XX}$ be a regular proper
excellent  two-dimensional scheme 
over the ring of integers in $k$ with the function field $K$. Let
$\alpha \in H^2(K, ~\mu_l)$ be a symbol  and $\zeta \in H^3(K, ~\mu_l)$. Suppose
that $\rama \cup \ramz$ is a union of regular curves with normal
crossings. If for every codimension one point $x$ of $\XX$ there
exists $f_x \in  K_x^*$ such that $\zeta = \alpha \cdot (f_x) \in
H^3(K_x, ~\mu_l)$,  then there exists $f \in K^*$ such that $\zeta =
\alpha \cdot (f)  \in H^3(K, ~\mu_l)$.
\end{cor}

\begin{proof} Let $D $ be an irreducible curve on ${\XX}$. Then the
residue field $\kappa(D)$ is either a local field or the function
field of a curve over a finite field.  Since ${\XX}$ is proper over
the ring of integers $k$,  by the class field theory it follows that
$H^2_{nr}(\kappa(D), \mu_l) = 0$.  The
residue field $\kappa(P)$ at   every closed point $P$ of ${\XX}$ is
a finite field and hence for every finite extension $\kappa_P$ of
$\kappa(P)$ we have $H^2(\kappa_P, \mu_l) = 0$.  Thus ${\XX}$ is
$l$-special and 
by (\ref{main}), there exists $f \in K^*$ such that $\zeta
- \alpha \cdot (f) \in H^3_{nr}(K/\XX, ~\mu_l)$. By a  theorem of
Kato (\cite{Ka}, 5.2), we have $H^3_{nr}(K/\XX, ~\mu_l) = 0$.  In particular
$\zeta = \alpha \cdot (f)$. \end{proof}

\begin{cor} 
\label{local-global-ps2}
  Let $k$ be a $p$-adic field
and $l$ a prime not equal to $p$. Let $X$ be a curve over $k$ and
$K$ its function field.     Let $\alpha \in H^2(K, ~\mu_l^{\otimes 2})$ be a
symbol and $\zeta \in H^3(K, ~\mu_l^{\otimes 3})$. 
  If for every discrete valuation $\nu$
of $K$ there
exists $f_\nu \in  K_\nu^*$ such that $\zeta = \alpha \cdot (f_\nu) \in
H^3(K_\nu, ~\mu_l^{\otimes 3})$,  then there exists $f \in K^*$ such that $\zeta =
\alpha \cdot (f)  \in H^3(K, ~\mu_l^{\otimes 3})$.
\end{cor}

\begin{proof}  Let $\rho$ be a primitive $l^{\rm th}$ root of
unity and $L = F(\rho)$. Suppose that
there exists $g \in L^*$ such that $\zeta = \alpha \cdot (g) \in H^3(L, \mu_l^{\otimes 3})$. 
Since $[L : K]$ is coprime to $l$, by taking corestrictions, we 
see that $\zeta = \alpha \cdot (f) \in H^3(K, \mu_l^{\otimes 3})$ for
some $f \in K^*$.    Let $\nu$ be a discrete valuation of 
$L$ and $\nu'$ be its restriction to $K$. Then, by the assumption, 
there exists $f_{\nu'} \in K^*$ such that $\zeta = \alpha \cdot (f_{\nu'})
 \in H^3(K_{\nu'}, \mu_l^{\otimes 3})$. Set $f_\nu = f_{\nu'} \in K_{\nu'} \subset L_\nu$. 
  $\zeta = \alpha \cdot (f_\nu) \in H^(L, \mu_l^{\otimes 3})$.
  Thus, replacing $K$ by $L$, we assume that $K$ contains a primitive
  $l^{\rm th}$ root of unity. 

Let ${\XX}$ be a regular proper
excellent  two-dimensional scheme 
over the ring of integers in $k$ with the function field $K$ such that
$\rama \cup \ramz$ is a union of regular curves with normal
crossings. For every  $x \in {\XX}^1$, by hypothesis, there exists
$f_x \in K_x$ such that $\zeta - \alpha \cdot (f_x)$ is unramified at
$x$.  The corollary follows by (\ref{local-global-ps1}). \end{proof}

\begin{cor} 
\label{local-global-ff1}
 Let $k$ be a finite field and $l$ a
prime not equal to the characteristic of $k$. Let $X$ be a smooth
projective surface over $k$ and $K$ its function field. Assume
that $K$ contains a primitive $l^{\rm th}$ root of unity.  Let
$\alpha \in H^2(K, ~\mu_l)$ be a symbol and $\zeta \in H^3(K, ~\mu_l)$.  Suppose
that $\rama \cup \ramz$ is a union of regular curves with normal
crossings. If for every codimension one point $x$ of $\XX$ there
exists  $f_x \in K_x^*$ such that $\zeta = \alpha \cdot (f_x) \in H^3(K_x,
~\mu_l)$, then there exists $f \in K^*$ such that $\zeta =
\alpha \cdot (f) \in H^3(K, ~\mu_l)$.
\end{cor}

\begin{proof}  For every irreducible curve $D$ on $X$,
$\kappa(D)$ is a function  field of curve over $k$ and for every
closed point $P$ of $X$, $\kappa(P)$ is a finite field. 
Thus, as in  the proof of (\ref{local-global-ps1}), $X$ is $l$-special.
By (\ref{main}), there exists $f \in K^*$ such that $\zeta
- \alpha \cdot (f) \in H^3_{nr}(K/X, ~\mu_l)$. By (\cite{CTSS} p.790, \cite{Ka}, Thm. 0.8),
we have $H^3_{nr}(K/X, ~\mu_l) = 0$.  In particular
$\zeta = \alpha \cdot (f)$. \end{proof}

\begin{cor} 
\label{local-global-ff2} Let $k$ be a finite field and $l$ a
prime not equal to the characteristic of $k$. Let $X$ be a smooth
projective surface over $k$ and $K$ its function field.     Let
$\alpha \in H^2(K, ~\mu_l^{\otimes 2})$ be a symbol and
 $\zeta  \in H^3(K, ~\mu_l^{\otimes 3})$.  
 If for every   discrete valuation $\nu$ of $K$  there
exists  $f_\nu \in K_\nu^*$ such that $\zeta = \alpha \cdot (f_\nu) \in H^3(K_\nu,
~\mu_l^{\otimes 3})$, then there exists $f \in K^*$ such that $\zeta =
\alpha \cdot (f) \in H^3(K, ~\mu_l^{\otimes 3})$.
\end{cor}

\begin{proof}   The proof follows on the same lines as the proof of  
(\ref{local-global-ps2}). \end{proof}

\begin{cor} 
\label{local-global-special}   Let $l$ be a prime and  ${\XX}$  an  $l$-special 
scheme.  Let $K$ be  the function field $\XX$. Let
$\alpha \in H^2(K, ~\mu_l)$ be a symbol  and $\zeta \in H^3(K,
~\mu_l)$. 
Suppose for every  point $z$ of $ {\XX}$, there exists $f_z \in K^*_z$
such that $\zeta  -  \alpha \cdot (f_z) \in H^3_{nr}(K_z/\hat{A}_z,  \mu_l)$.
Then  there exists $f \in K^*$ such that $\zeta -  \alpha \cdot (f)
\in H^3_{nr}(K/{\XX}, \mu_l)$. 
\end{cor}

\begin{proof}  Let $\eta : {\XX}' \to {\XX}$ be a blow-up. 
Let $x \in \XX'$ be a codimenson one point and $z = \eta(x) \in \XX$. 
By the hypothesis there exists $f_z \in K_z^*$ such that 
$\zeta -  \alpha \cdot (f_z) \in H^3_{nr}(K_z/\hat{A}_z, \mu_l)$.

Suppose that $z $ is a codimension one point of $\XX$.
Then $\OO_{\XX', x} \simeq \OO_{\XX, z}$ and hence
$\zeta - \alpha \cdot (f_z)$ is unramified at $x$. 
Suppose that $z$ is  a codimension two point of $\XX$.
Then, by (\ref{purity}), $\zeta - \alpha \cdot (f_z)$ is unramified at $x$.

Hence by replacing ${\XX}$ by a blow-up,
we assume that the support of $\rama \cup \ramz$ is a union of regular
curves with normal crossings. 
 By ( \ref{main}), there exists $f \in K^*$ such that $\zeta - \alpha \cdot
(f) \in H^3_{nr}(K/{\XX}, \mu_l)$.    \end{proof}

\section{The $H^3$ of the function field of a surface}
 
Let ${\XX}$ be an excellent  regular
integral surface  and   $K$
its function field. Let  $l$ be a prime which is a
unit on ${\XX}$.  Suppose that $K$ contains a primitive $l^{th}$ root
of unity.   In this section,  under certain conditions on $K$, we show that
every element in $H^3(K, \mu_l)$ is a symbol up to an 
element in  $H^3_{nr}(K/\XX, \mu_l)$. We begin with the
following.

\begin{lemma} 
\label{symbol-local}
Let $A$ be  a complete regular local ring of
dimension two with the field of fractions $K$ and the residue field
$\kappa$.  Let $l$ be a prime not equal to char$(\kappa)$. Suppose that
$K$ contains a primitive $l^{\rm th}$ root of unity 
 and $H^2(\kappa, \mu_l) = 0$.
Let $m = (\pi, \delta)$
be  the  maximal ideal $m $ of $A$.  Let $\zeta \in
H^3(F, \mu_l)$ be such that $\zeta$ is ramified on $A$ at most at
$\pi$ or $\delta$.  Let $u \in A$ be a unit such that
$\partial_{\overline{\delta}}(\partial_{\pi}(\zeta)) = (\tilde{u}) \in
H^1(k, \mu_l)$, where $\bar{~}$ denotes  the image modulo $\pi$ and
$\tilde{~}$ denotes  the image modulo the maximal ideal.  
Then $\zeta - (u) \cdot (\delta) \cdot (\pi)$ is unramified on $A$.
\end{lemma}

\begin{proof}   Since $\partial_{\overline{\delta}}
(\partial_{\pi}(\zeta)) = (\tilde{u})$ and $A/(\pi)$ is a complete
discrete valuation ring with $\overline{\delta}$ as a parameter, 
 $\partial_{\pi}(\zeta) - (\overline{u}) \cdot
(\overline{\delta})$ is unramified at $\overline{\delta}$.  Since
$H^2(\kappa, \mu_l) = 0$, $\partial_{\pi}(\zeta) = (\overline{u})
\cdot (\overline{\delta})$. 
Since $\partial_{\pi}((u) \cdot (\delta) \cdot
(\pi)) = (\overline{u}) \cdot (\overline{\delta})$, $\zeta - (u) \cdot
(\delta) \cdot (\pi)$ is unramified except possibly  at $\delta$. Hence, by
(\ref{ramified-at-pi}), $\zeta - (u) \cdot (\delta) \cdot (\pi)$ is unramified on
$A$.    
\end{proof}

\begin{theorem}
\label{symbol-ur}
 Let $l$ be  a prime and  ${\XX}$ an
$l$-special scheme.  Let $K$ be the  function field   of ${\XX}$.   
Suppose that $K$ contains a primitive $l^{th}$ root
of unity.    Then every element $\zeta
\in  H^3(K, ~\mu_l)$ is of the form $\zeta' + (a) \cdot (b) \cdot (c)$
for some $\zeta' \in H^3_{nr}(K/\XX, \mu_l)$.  
\end{theorem}

\begin{proof} We note that for any blow-up ${\XX}' \to {\XX}$,
$H^3_{nr}(K/{\XX}',  \mu_l) \subset H^3_{nr}(K/{\XX}, \mu_l)$, 
enabling us to replace $\XX$ by any blow-up.  
Let $\zeta \in H^3(K, ~\mu_l)$.  By replacing
${\XX}$ by  a blow-up, we assume that $\ramz  = C \cup E$ with $C$
and $E$ regular curves (not necessarily irreducible)  on ${\XX}$
with normal crossings.  Suppose $C = \cup_i C_i$ and $E = \cup_j E_j$
with $C_i$ and $E_j$ regular irreducible curves on ${\XX}$.
Let ${\PP}$ be the  finite set of closed points of $C \cap E$.  Let
$B$ be the semi-local ring at ${\PP}$, $\pi_i \in B$ be a prime
defining $C_i$ on $B$ and $\delta_j \in B$ a prime defining $E_j$ on
$B$.    For each $P \in {\PP} \cap C_i$, let $u_P \in
\kappa(P)^*$ be such that $\partial_P(\partial_{C_i}(\zeta)) = (u_P)$. Let $u_i
\in B/(\pi_i)$ be such that $ u_i(P) =  u_P \in \kappa(P)$  
for all $P \in C_i \cap {\PP}$.

Let $\nu$ be a discrete valuation of $\kappa(C_i)(\sqrt[l]{u_i})$ which
is centered on a closed point $P$ of $C_i$. 
Suppose that $P \not\in  {\PP}$.  Then, by (\ref{ramified-at-pi}), $\zeta$
is unramified on $\hat{A}_P$ and hence, by (\ref{residue-unramified-at-P}), 
$\partial_{C_i}(\zeta) \in
H^2(\kappa(C_i), \mu_l)$ is unramified at $P$. In particular, 
 $\partial_{C_i}(\zeta) \otimes \kappa(C_i)(\sqrt[l]{u})$
is unramified at  $\nu$. Suppose that $P \in \PP$.  
Then, by the choice of $u_i$, we have $\partial_P(\partial_{C_i}(\zeta))
= (u_P) = (u_i(P))$.  Since $\partial_\nu (\partial_{C_i}(\zeta 
\otimes \kappa(C_i)(\sqrt[l]{u_i}))) =  \partial_P(\partial_{C_i}(\zeta)) \otimes
\kappa(\nu) = (u_i(P)) \otimes \kappa(\nu)$ and $u_i(P) $ is 
an $l^{\rm th}$ power in $\kappa(\nu)$, $\partial_\nu(\partial_{C_i}(\zeta)
\otimes  \kappa(C_i)(\sqrt[l]{u_i}) )$
is trivial.
Hence $\partial_{C_i}(\zeta) \otimes  \kappa(C_i)(\sqrt[l]{u_i})$ 
is unramified at every discrete valuation of
$\kappa(C_i)(\sqrt[l]{u_i})$ centered on a closed point of $C_i$. 
Since  $\XX$ is $l$-special,  $H^2_{nr}(\kappa(C_i)(\sqrt[l]{u_i})/C_i, \mu_l) = 0$ and
$\partial_{C_i}(\zeta) $ is zero over $\kappa(C_i)(\sqrt[l]{u_i})$.
Hence, $\partial_{C_i}(\zeta) = (u_i)
\cdot (a_i)$ for some $a_i \in \kappa(C_i)$.   Similarly, let $v_j \in
B/(\delta_j)$ be a unit  such that $(v_j (P))
= \partial_{P}(\partial_{E_j}(\zeta))$ for all $P \in E_j \cap {\PP
  }$ and $\partial_{E_j}(\zeta) = (v_j) \cdot (b_j) $
for some $b_j \in \kappa(E_j)$  and $v_j \in B/(\delta_j)$. 

Let $P \in C_i \cap E_j$ for some $i$ and $ j$. Since $
\partial_P(\partial_{C_i}(\zeta))
+ \partial_P(\partial_{E_j}(\zeta)) = 0 $ (cf. (\ref{complex})),  $u_i(P) =
v_j(P)^{-1}$. By (\cite{S3},
0.3(b)), there exists  $u \in B$ such that $u = u_i$ modulo $\pi_i$
for all $i$ and $u = v_j^{-1}$  modulo $\delta_j$ for all $j$ .  Let $b =
\prod_i \pi_i  \prod_j \delta_j\in B$ and $\alpha = (u) \cdot (b)$. 
Then for $P \in C_i \cap E_j$, we have $b = w_P\pi_i\delta_j$
for unit $w_P$ at $P$. 
 We now show that for every codimension one point 
$D$ of $\XX$,  there exists $f_D  \in
K_{D}^*$ such that $\zeta - \alpha \cdot (f_D )$ is unramified at
$D$.  
  
Let $D$ be a codimension one point of 
${\XX}$. If $D \neq C_i$ for all $i$ and $D \neq E_j$ for all $j$,
then $\zeta$ is unramified at  $D$ and hence $\zeta  = 0 = \alpha
\cdot (1)$ over $K_D$.  Suppose $D = C_i$ for some $i$.  We have
$\partial_{C_i}(\zeta) = (u_i) \cdot  (a_i)$ for some $a_i \in
\kappa(C_i)^*$.  Let $f_i \in B$ be a unit at $C_i$ such that $f_i$ maps to
$a_i$ in $\kappa(C_i)$.  Since $\partial_{C_i}(\alpha \cdot (f_i^{-1})) = 
\partial_{C_i}((u) \cdot (b) \cdot (f_i^{-1}) = 
(u_i) \cdot (a_i)$, we have $\zeta - \alpha \cdot (f_i^{-1})$ is
unramified at $C_i$.  For  $D = E_j$,   arguing in a similar way, one
finds   $g_j \in K$  such that   $\zeta - \alpha \cdot (g_j)$ is
unramified at $E_j$.

Let $\eta : \XX' \to \XX$ be a blow-up such that $\ramap \cup
\ramzp$ is a union of regular curves with normal crossings. Let $x$ be
a codimension one point of $\XX'$.  We show that there exists $f_x
\in K_x^*$ such that $\zeta - \alpha \cdot (f_x) \in H^3_{nr}(K_x,
\mu_l)$.  If $\eta(x) = y$ is a codimension one point of $ \XX$,
then $K_y = K_x$ and hence there exists $f_x \in K_x^*$ such that
$\zeta - \alpha \cdot (f_x) \in H^3_{nr}(K_x, \mu_l)$. Let $\eta(x) =
P$ be a closed point of $\XX$.  If $P \not\in C_i \cap E_j$ for all
$i $ and $j$, then   by (\ref{ramified-at-pi}) 
$\zeta$ is unramified  on $\hat{A}_P$.
One concludes using (\ref{purity}) that   $\zeta$ is unramified 
over $K_x$. In particular $\zeta - \alpha
\cdot (1)$ is unramified at $x$. Suppose $P \in C_i \cap E_j$ for some
$i $ and $ j$.  Since $\partial_P(\partial_{C_i}(\zeta)) = (u(P))$, by
(\ref{symbol-local}), we have $\zeta - (u) \cdot (\delta_j) \cdot (\pi_i)$ is
unramified on $\hat{A}_P$.   We have 
 
$$
\begin{array}{rcl}
\zeta - \alpha \cdot (\pi_i) & = &\zeta -  (u) \cdot (b) \cdot (\pi_i) \\
& = &  \zeta - (u) \cdot (w_P\pi_i\delta_j) \cdot (\pi_i) \\
& = & \zeta - (u) \cdot (-w_P\delta_j) \cdot (\pi_i)  \\
& = & \zeta - (u) \cdot (\delta_j) \cdot (\pi_i) - (u) \cdot (-w_P) \cdot (\pi_i) 
\end{array}
$$
Since $(u) \cdot (-w_P) = 0$ in $H^2(K_P, \mu_l)$, 
$\zeta - \alpha \cdot (\pi_i) = \zeta - (u) \cdot (\delta_j) \cdot (\pi_i)$
is  unramified on $\hat{A}_P$.  Once again by (\ref{purity}), we conclude that
$\zeta - \alpha \cdot (\pi_i)$ is unramified at $x$. 

Thus, by (\ref{main}),  there exists $a
\in K^*$ such that $\zeta - \alpha \cdot (a) = \zeta- (u) \cdot (b) \cdot
(a) \in H^3_{nr}(K/\XX, \mu_l)$.  \end{proof}

\begin{cor} 
\label{symbol-p-adic}
(cf. \cite{PS2}, 3.5) Let $k$ be a $p$-adic field and
$l$ a prime not equal to $p$. Suppose that $k$ contains a primitive
$l^{\rm th}$ root of unity. Let $X$ be a curve over $k$ and $K$ its
function field. Then every element in $H^3(K, ~\mu_l)$ is a symbol.
\end{cor}

\begin{proof} Let ${\OO}_k$ be the ring of integers in $k$
and ${\XX}$  a regular integral surface, proper over Spec$({\OO
}_k)$ with function field $K$.   We have seen
(in the proof of  (\ref{local-global-ps1})), that $\XX$  is $l$-special. 
Since $H^3_{nr}(K/\XX, ~\mu_l) = 0$ (\cite{Ka}, 5.2), the corollary follows
from (\ref{symbol-ur}). \end{proof}

\begin{cor} 
\label{symbol-ff}
 Let $k$ be a finite  field and $X$ a smooth,
projective, geometrically integral surface over $k$.  Let $K$ be the
function field of $X$. Let $l$ be a prime not equal to char$(k)$.
Assume that $k$ contains a primitive $l^{th}$ root of unity. Then every
element in $H^3(K, ~\mu_l)$ is a symbol.
\end{cor}

\begin{proof}    We have seen (  in the proof of (\ref{local-global-ff1})) that
 $X$ is $l$-special.   Since
$H^3_{nr}(K/\XX, ~\mu_l)  \\ = 0$ (\cite{Ka} Thm. 0.8), the corollary follows
from (\ref{symbol-ur}). \end{proof}

\section{ The unramified $H^3$ of  function fields of  conic
fibrations over surfaces}

Let $K$ be a field of characteristic not equal to 2.
We say that a {\it local-global principle holds  for $H^3(K, ~\mu_2)$ 
in terms of symbols in $H^2(K, ~\mu_2)$} if given  $\zeta \in
H^3(K, ~\mu_2)$, $\alpha \in H^2(K, ~\mu_2)$ representing a quaternion
algebra over $K$, if  for every discrete valuation $\nu$ of $K$, there
exists $f_{\nu} \in K_{\nu}^*$ such that $\zeta - \alpha \cdot
(f_{\nu}) \in H^3_{nr} (K_{\nu}, ~\mu_2)$  then there exists $f \in
K^*$ such that $\zeta - \alpha \cdot (f) \in H^3_{nr}(K, ~\mu_2)$.

\begin{lemma}
\label{dvr-conic} 
Let $K$ be a field of characteristic not equal to 2
and $\alpha$ a quaternion algebra over $K$.
Let  $C$ be the smooth projective conic associated to $\alpha$. 
Let  $\nu$ be  a   discrete valuation on $K$ with residue field $\kappa$.
Then there exists a discrete valuation $w$ on the function field $K(C)$
of $C$ with residue field $\kappa(w)$ satisfying: \\
1) If $\alpha$ is unramified at $\nu$,   then    $\kappa(w)$ is isomorphic to
the function field of the conic $\overline{C}$ associated to the specialisation
$\overline{\alpha}$ of $\alpha$. \\
2) If $\alpha$ is ramified at $\nu$, then $\kappa(w)$ is isomorphic to 
the rational function field in one variable over the quadratic extension of $\kappa$ 
given by the residue of $\alpha$.
\end{lemma}

\begin{proof} 
 Let $\alpha  = (a) \cdot (b)$ for some units $a, b \in K^*$.
Then  $K(C) = K(x)[y]/(y^2 - ax^2 - b)$.  
Let $R$ be the   valuation ring of $\nu$ and $\pi$  a parameter in $R$.
Let $\tilde{\nu}$ be the discrete valuation of $K(x)$ given by 
the height one prime ideal $(\pi)$ of $R[x]$.
Then the residue field of $\tilde{\nu}$ is $\kappa(x)$.

 Suppose that $\alpha$ is unramified at $\nu$.
Then $ a$ and $ b$ may  be chosen to be  units  in  $R$ (cf. \ref{dvr-symbol})
and $\overline{\alpha} = (\overline{a}) \cdot (\overline{b})$.
Since $a$ and $b$ are units of $R$, the quadratic extension 
 $K(C)$ of $K(x)$ is unramified at $\tilde{\nu}$. 
 Hence there is a unique extension $w$ of $\tilde{\nu}$ to $K(C)$
  with residue field $\kappa(x)[y]/(y^2 - \overline{a}x^2 - \overline{b})$.

 Suppose that $\alpha$ is ramified at $\nu$. We can  choose $a$  to be 
 a unit and $b$ to be a paramter in $R$  (cf. \ref{dvr-symbol}); the residue of
  $\alpha$ is  $(\overline{a})$.  Since $\alpha$ is ramified at $\nu$,
  $\overline{a}$ is not a square in $\kappa(\nu)$. Hence there is a unique extension 
  $w$ of $\tilde{\nu}$ to $K(C)$  with residue field $\kappa(x)[y]/(y^2 - \overline{a})$.
  \end{proof}

\begin{theorem}
\label{lg-ur}
 Let $K$ be a field of characteristic
 not equal to 2. Suppose that
for every discrete valuation $\nu$ of $K$, the characteristic of
$\kappa(\nu)$ is not equal to $2$. The local-global principle for
$H^3(K, ~\mu_2)$ in terms of symbols in $H^2(K, ~\mu_2)$ holds if
and only if for every smooth conic $C$ over $K$, the restriction map
$H^3_{nr}(K, ~\mu_2) \to (H^3_{nr}(K(C), ~\mu_2) \cap {\rm
image}(H^3(K, \mu_2)))  $ is onto.
\end{theorem}

\begin{proof} Suppose that the local-global principle for
$H^3(K, ~\mu_2)$ in terms of symbols in $H^2(K, ~\mu_2)$ holds. Let
$\beta \in H^3_{nr}(K(C), ~\mu_2)$ which is in the image of $H^3(K,
\mu_2)$. Let $\zeta \in H^3(K, ~\mu_2)$ map to $\beta$. Let
$\alpha \in H^2(K, ~\mu_2)$ be the element representing the
quaternion algebra associated to the conic $C$.

 For any $f \in K^*$,   $\alpha \cdot (f) \in H^3(K, \mu_2)$ maps to zero in
$H^3(K(C), \mu_2)$.  Thus it is enough to show that
there exists $f \in K^*$ such that $\zeta - \alpha \cdot (f) \in
H^3_{nr}(K, ~\mu_2)$. Since the local-global principle for $H^3(K,
~\mu_2)$ in terms of symbols in $H^2(K, ~\mu_2)$ holds, it suffices
to show that for every discrete valuation $\nu$ of $K$, there exists
$f_{\nu} \in K_{\nu}^*$ such that $\zeta - \alpha \cdot (f_{\nu})
\in H^3_{nr}(K_{\nu}, ~\mu_2)$.

Let $\nu$ be a  discrete valuation of $K$.   If $\zeta$ is unramified
at ${\nu}$, then $\zeta - \alpha \cdot (1) \in H^3_{nr}(K_{\nu}, \mu_2)$.
Assume that $\zeta$ is ramified at ${\nu}$; i.e. $\partial_{\nu}(\zeta) \neq 0$.

Suppose that $\alpha$ is unramified at $\nu$. Then, by (\ref{dvr-conic}), 
there exists a
discrete valuation $w$ on $K(C)$ extending $\nu$ on $K$ with residue
field $\kappa({\nu})(\overline{C})$, where $\overline{C}$ is the
reduction of $C$ at $\nu$. Since $\beta \in H^3_{nr}(K(C), ~\mu_2)$
is the image of $\zeta \in H^3(K, ~\mu_2)$, we have
$\partial_{\nu}(\zeta) = \partial_w(\beta) = 0$ over
$\kappa(\nu)(\overline{C})$.   Since $\partial_{\nu}(\zeta) \neq 0$
in $H^2(\kappa(\nu), ~\mu_2)$ and $\overline{C}$ is the conic
associated to the specialisation $\overline{\alpha}$ of $\alpha$ at
$\nu$, by a theorem of Amitsur (\cite{amitsur}, Theorem 9.3), we have $
\partial_{\nu}(\zeta) =
\overline{\alpha}$. Since $K_{\nu}$ is a complete  discretely   valued
field with residue field $\kappa(\nu)$, we have  $\zeta - \alpha
\cdot (\pi) \in H^3_{nr}(K_{\nu}, ~\mu_2)$.

Suppose that $\alpha$ is ramified at $\nu$. Then  
$\alpha = (a) \cdot (b\pi)$ over $K_{\nu}$ for some $a, b \in
K_{\nu}^*$ which are units in ${\OO}_{\nu}$ (cf. \ref{dvr-symbol}). 
 By (\ref{dvr-conic}), 
there is a discrete valuation  $w$ on $K(C)$ extending   $\nu$  
 with residue field
$\kappa(\nu)(\sqrt{\overline{a}})(t)$. As above, since the residue
$\partial_{w}(\zeta)$  is zero over $\kappa(w) =
\kappa(\nu)(\sqrt{\overline{a}})(t)$, it follows that
$\partial_{\nu}(\zeta)$ is split over $K(\nu)(\sqrt{\overline{a}})$
and hence $\partial_{\nu}(\zeta) = (\overline{a}) \cdot
(\overline{c})$ for some unit $c \in {\OO}_{\nu}$. Hence $\zeta'
= \zeta - (a) \cdot (c) \cdot (\pi)$ is unramified at $\nu$. We have
over $K_{\nu}$
$$
\zeta - \alpha \cdot (c)  =    \zeta - (a) \cdot (c) \cdot  (\pi) +
(a) \cdot (c) \cdot (\pi) - (a) \cdot (c) \cdot (b\pi)   = \zeta' -
(a) \cdot (c) \cdot (b).
$$
Since $a, b, c$ are units at $\nu$,  $(a) \cdot (c) \cdot (b)$ is
unramified at $\nu$. In particular $\zeta - \alpha \cdot (c)$ is
unramified at $\nu$.

Conversely, suppose that the restriction map $ H^3_{nr}(K, ~\mu_2)
\to H^3_{nr}(K(C), ~\mu_2) \\ \cap {\rm image}(H^3(K, \mu_2))$ is
onto.

Let $\zeta \in H^3(K, ~\mu_2)$ and $\alpha \in H^2(K, ~\mu_2)$
representing a quaternion algebra. Suppose that for every discrete
valuation $\nu$ there exists $f_{\nu} \in K_{\nu}^*$ such that
$\zeta - \alpha \cdot (f_{\nu}) \in H^3_{nr}(K_{\nu}, \mu_2)$.

Let $C$ be the conic associated to $\alpha$. Let $w$ be a discrete
valuation on $K(C)$. Suppose that $w$ is trivial on $K$. Since
$\zeta $ is defined over $K$, $\zeta$ is unramified at $w$. 
 Therefore we, assume
that  the valuation $w$ is non-trivial when restricted to  $K$. Let $\nu$ be the restriction  of
$w$ to $K$. Then we have $K_{\nu} \subset K_{\nu}(C) \subset
K(C)_w$. Since $\zeta - \alpha \cdot (f_{\nu})$ is unramified at
$\nu$ and $\alpha \cdot (f_{\nu})$ is zero over $K(C)_w$, $\zeta$ is
unramified at $w$.

In particular $\zeta \in H^3_{nr}(K(C), ~\mu_2) \cap {\rm
image}(H^3(K, \mu_2))$. By the hypothesis, there exists $\beta \in
H^3_{nr}(K,~\mu_2)$ which maps to $\zeta$ in $H^3(K(C), ~\mu_2)$.
Since $\zeta - \beta \in H^3(K, ~\mu_2)$ maps to $0$ in $H^3(K(C),
~\mu_2)$, we have $\zeta = \beta + \alpha \cdot (f)$ for some $f \in
K^*$ (\cite{arason}, Corollar 5.5). Since $\beta \in H^3_{nr}(K, ~\mu_2)$ we are done.
\end{proof}

\vskip 5mm

{\it We thank A. Pirutka for the observation that no condition
on the Brauer group of the residue fields is required in the
hypothesis of (\ref{lg-ur}).}

\begin{theorem}
\label{h3ur-ff}
 Let $k$ be a finite  field of characteristic
not equal to 2. Let $X$ be a smooth, projective, integral surface
over $k$ and $K$ its function field. Let $n$ be a natural number
coprime to the characteristic of $K$ and $C$ a smooth conic over $K$. Then
$H^3_{nr}(K(C)/k, ~\mu_n^{\otimes 2}) = 0$.
\end{theorem}

\begin{proof}   By  (\cite{CTSS},  p.790 and \cite{Ka}, 
Thm. 0.8), we have    $H^3_{nr}(K, \mu_n^{\otimes 2}) = 0$ and $H^3_{nr}(K, {\bf Q}/{\bf Z}(2)) = 0$.
We shall use this fact without further reference. 
 
Suppose that $C$ is the projective line.  
Since  $H^3_{nr}(K(C)/C, \mu_n^{\otimes 2}) = H^3(K, \mu_n^{\otimes 2})$ 
(cf. \cite{CT1}, Prop. 4.1.4)  
and every discrete valuation of $K$ extends to a discrete valuation of  $K(C)$,    
$H^3_{nr}(K(C), ~\mu_n^{\otimes 2})  \simeq H^3_{nr}(K, \mu_n^{\otimes 2})$. Hence 
$H^3_{nr}(K(C), \mu_n^{\otimes 2}) =  0$.

Suppose  that $C$ is an anisotropic conic. Since $C$ is isomorphic to
the projective line over a separable quadratic  extension of $K$,  
we conclude that $H^3_{nr}(K(C), ~\mu_n^{\otimes 2})$ is 2-torsion. Hence $H^3_{nr}(K(C),
~\mu_n^{\otimes 2}) \subset H^3_{nr}(K(C), ~\mu_2^{\otimes 2})$, both being regarded as subgroups of
$H^3(K(C), {\bf Q}/{\bf Z}(2))$ (\cite{MS}).  Thus it is enough to show
that  $H^3_{nr}(K(C), ~\mu_2) = 0$.

Let $\beta \in H^3_{nr}(K(C), ~\mu_2)$. 
Since $C$ is a smooth conic, there exists $\zeta \in H^3(K, {\bf
Q}/{\bf Z}(2))$ which maps to the image of $\beta$ in $H^3(K(C), {\bf
Q}/{\bf Z}(2))$  (\cite{Su}). Let $v$ be a
discrete valuation on $K$.  By, (\ref{dvr-conic}), there exists a
discrete valuation $w$ of $K(C)$ with $\kappa(w)$ either the
function field of a conic over $\kappa(v)$ or the function field of
the projective line over a quadratic extension of $\kappa(v)$.
Since $\partial_v(\zeta) \otimes \kappa(w) =
\partial_w(\beta) = 0$,  it follows that
$2\partial_v(\zeta) = 0$. Hence $2\zeta$ is unramified at $v$. Since
$H^3_{nr}(K, ~{\bf Q}/{\bf Z}(2)) = 0$, $2\zeta = 0$,
i.e. $\zeta \in H^3(K, \mu_2)$. Hence $\beta \in H^3_{nr}(K(C),
\mu_2) \cap {\rm image}(H^3(K, \mu_2))$. By (\ref{main}) and (\ref{lg-ur}), the map
$H^3_{nr}(K, ~\mu_2) \to H^3_{nr}(K(C), ~\mu_2) \cap {\rm
image}(H^3(K, \mu_2))$ is onto. Since $H^3_{nr}(K,~\mu_2) = 0$,  $\beta = 0$ and 
the theorem follows. \end{proof}

\begin{theorem} 
\label{h3nr-p-adic}
 Let $k$ be a $p$-adic field with $p \neq
2$. Let $X$ be a smooth, projective, integral curve over $k$ and $K$
its function field. Let $l$ be a prime different from  2 and $p$,  and $C$ a smooth 
conic over $K$. Then $H^3_{nr}(K(C), ~\mu_l) = 0$.
\end{theorem}

\begin{proof} Since $H^3_{nr}(K,~\mu_l) = 0$ (\cite{Ka}, 5.2),
the theorem follows from (\ref{main}) and (\ref{lg-ur})
 following the proof of (\ref{h3ur-ff}) . \end{proof}

\section{ Local-global principle for zero-cycles and integral
Tate conjecture}

Let $k$ be a global field and $X$ a smooth, projective, geometrically 
integral variety over $k$. It is a conjecture 
due to Colliot-Th\'el\`ene--Sansuc and Kato-Saito (cf. \cite{CT2}, \S 4) 
that 
the Brauer-Manin obstruction is the only obstruction to the local-global
principle  for the existence of zero-cycles of degree one on $X$.
We show that for a certain class  of surfaces over a global field
of  positive characteristic, this conjecture is true. 
A corresponding result for conic fibrations of the projective line
over a number field is due to Salberger (\cite{Sb}, Corollary on p.505).

For more details on the results of this section we refer to   the  paper of  
Colliot-Th\'el\`ene and Bruno Kahn (\cite{CTK}). 
  
Let  ${ {\F}}$ be a finite field and 
$C$  a smooth, projective, geometrically integral curve over
${ {\F}}$.  Let $X$ be a smooth projective, geometrically
integral variety over ${ {\F}}$ of dimension $d + 1$ with a
flat morphism $X \to C$ with generic fibre
$X_{\eta}/{ {\F}}(C)$  smooth and geometrically integral. Let
$l$ be a prime not equal to the characteristic of ${ {\F}}$. A
theorem of Saito   (\cite{Sa},  see also \cite{CT2}, Proposition 3.2) 
asserts that if
the cycle map $CH^d(X) \otimes { {\Z}}_l \to H^{2d}(X,
{ {\Z}}_l(d))$ is onto, then, the Brauer-Manin obstruction  is
the only obstruction to the local-global principle for the
existence of zero-cycles of degree coprime to $l$ on $X_{\eta}$; 
more precisely, if there exists a family of local 0-cycles of degree prime to $l$, 
orthogonal to the Brauer group of $X_{\eta}$, then there exists on 
$X_{\eta}$ a 0-cycle of degree prime to $l$.

For a certain class $B_{Tate}({ {\F}})$ of smooth, projective
varieties $Y$, Bruno Kahn (\cite{K}) showed that the cycle map $CH^2(Y) \otimes
{ {\Z}}_l \to H^4(Y,{ {\Z}}_l(2))$ is onto if and only if
$H^3_{nr}({ {\F}}(Y)/\F, { {\Q}}_l/{ {\Z}}_l(2)) = 0$.
Suppose $Y$ is a smooth projective variety with a dominant morphism
$Y \to X$, where $X$ is a smooth geometrically ruled surface, with
the generic fibre a smooth conic over ${ {\F}}(X)$. By a
theorem of Soul\'e (\cite{So}, Example 3.3) $Y$ is in $B_{Tate}({ {\F}})$. 
We have the following

\begin{theorem}
\label{cyclemap}
 Let ${ {\F}}$ be a field of characteristic not
equal to 2. Let $l$ be a prime not equal to the characteristic of
${ {\F}}$. Let $X$ be a smooth, projective, geometrically ruled
surface over ${ {\F}}$ and $Y \to X$ a surjective morphism with
the generic fibre a smooth conic over ${ {\F}}(X)$. Then the
cycle map $CH^2(Y) \otimes { {\Z}}_l \to H^4(Y,
{ {\Z}}_l(2))$ is surjective.
\end{theorem}

\begin{proof} Since $Y$ is in $B_{Tate}(\F)$, 
the result follows from a combination of Kahn's theorem
 (\cite{K}, Corollaire 3.13) and 
the well-known exact sequences
$$H^3_{nr}(\F(Y)/\F,\mu_{l^n}^{\otimes 2}) \to CH^2(Y)/l^n) \to H^4_{et}(Y,\mu_{l^n}^{\otimes 2})$$
(cf.  \cite{CTK}, Section 3.2, page 17).
 \end{proof}

\begin{cor} (cf. \cite{CTK}, Corollary 7.9)
\label{bmob}
 Let ${ {\F}}$ be a finite field of
characteristic not equal to 2. Let $X$ be a smooth, projective,
geometrically ruled surface over ${ {\F}}$ and $Y$ be a smooth,
geometrically integral 3-fold over ${{\F}}$ with a surjective
morphism $Y \to X$ with generic fibre a smooth conic over
${ {\F}}(X)$. Let $C$ be a smooth, projective geometrically
integral curve over ${ {\F}}$ together with a surjective
morphism $X \to C$. Suppose that the generic fibre $X_{\eta}$
of $X \to C$ has a 0-cycle of degree a power of 2.
Let $Y_{\eta}$ be the generic fibre of the
composite morphism $Y \to X \to C$. Then the Brauer-Manin
obstruction is the only obstruction to the local global principle for
the existence of zero-cycles of degree one on $Y_{\eta}$.
\end{cor}

\begin{proof} Suppose that there is no Brauer-Manin
obstruction for the existence of zero-cycles of degree one on
$Y_{\eta}$. Suppose that $Y_{\eta}$ has a zero-cycle of degree one
over each completion of ${ {\F}}(C)$. Since the
characteristic of ${ {\F}}$ is not equal to 2, by (\ref{cyclemap}) and
(\cite{CT2} Proposition 3.2, \cite{Sa}), there exists a zero-cycle on $Y_{\eta}$ of odd degree.
Since $X_{\eta}$ has a 0-cycle of degree a power of 2 and 
 $Y_{\eta} \to X_{\eta}$ is a conic fibration, 
$Y_{\eta}$ has a 0-cycle of degree a power of  2. Hence
$Y_{\eta}$ has a zero-cycle of degree 1. \end{proof}
 
\begin{example} Let $E$ be a smooth, projective
geometrically integral non-rational curve over a finite field
${ {\F}}$ of characteristic not equal to 2. Let $X = E \times
{ {\P}}^1$.  Let $Y \to X$ be as in (\ref{bmob}) and $X \to
{ {\P}^1}$ be the projection. Let $Y_{\eta}$ be the generic
fibre of the composite $Y \to X \to { {\P}}^1$. Then $Y_{\eta}$
is a non-rational surface over ${ {\F}}(t)$ and the Brauer-Manin
obstruction is the only obstruction to local global principle for
the existence of zero-cycles of degree one on $Y_{\eta}$.
\end{example}

\section{ Appendix}
\noindent
{\bf Erratum to ``The $u$-invariant of the function fields of
  $p$-adic curves'' (Annals of Mathematics, 172 (2010), 1391-1405) by R. Parimala and V. Suresh }

\vskip 5mm

\noindent
 A hypothesis is missing in the statement of  (\cite{PS2}, Proposition 3.2 and Theorem 3.4).
This is bridged in (4.5) of this paper and using it  we complete the proof of $u(K) = 8$
 for function fields of non-dyadic $p$-adic curves (\cite{PS2}, 4.6).   We begin with the following 
  
\begin{lemma}
\label{app-lemma1}
Let $A$ be a complete regular local ring of dimension 2 with
the maximal ideal $m = (\pi, \delta)$,  the field of fractions 
$K$ and the residue field $\kappa$. 
Suppose that 2 invertible in $A$.
Let $q = \phi_0 \perp \phi_1\pi \perp \phi_2\delta
\perp \phi_3\pi\delta$ with $\phi_i$, $0 \leq i \leq 3$, regular quadratic forms over $A$. 
Let $g \in A$ be such that $ g$ is either a unit or $u\pi$ modulo $\delta$
and $v\delta$ modulo $\pi$  for some units 
$u, v \in A$. If $q$ represents $ g$  over   $K_\pi$ and $K_\delta$, then
$q$ represents $g$  over $K$. 
 \end{lemma}

\begin{proof}  We assume without loss of generality 
$q$ is anisotropic over $K$. Since $A$ is complete, the forms
$\phi_i$, $0 \leq i \leq 3$, are anisotropic modulo $m$. 
Suppose $g$ is a unit modulo $\pi$ or $\delta$. Then 
$g$ is a unit in $A$.  We have 
$q \perp <-g>  = \phi'_0 \perp   \phi_1\pi \perp \phi_2\delta
\perp \phi_3 \pi\delta$ with $\phi'_0 = \phi_0 \perp <-g> $
regular on $A$.  By taking residues at $\pi$ or $\delta$ and 
using the fact that  $A$ is complete  and 
$q \perp <-g>$ is isotropic over $K_\pi$ and
$K_\delta$, we conclude that $\phi'_0$  is isotropic over $K$ 
and hence $q \perp <-g>$ is isotropic over $K$ and  $q$ represents $g$ over $K$.

Assume that $g = u\delta$ modulo $\pi$ and
$g = v\delta$ modulo $\pi$.  Then $g = u\pi + v\delta + \pi\delta \theta$ 
for some $\theta \in A$.  Thus $g = (u + \delta\theta)\pi + v\delta$
and hence, by replacing $u$ with $u + \delta\theta$, we assume that
$g = u\pi + v\delta$ with $u, v$ units in $A$. Since 
$$q \perp <-g> =  \phi_0 \perp \phi_1\pi \perp \phi_2\delta
\perp \phi_3 \pi\delta - <u\pi + v\delta> $$
is isotropic over $K_\pi$ and $K_\delta$, 
by taking the residues at $\pi$ and $\delta$, we see that
$\phi_1$ represents $u$ and $\phi_2$ represents $v$
and hence $q$ represents $g = u\pi + v\delta$.
\end{proof}
 
\begin{lemma}
\label{app-lemma2}
 Let $A$ be a complete regular local ring of
dimension two with maximal ideal $m = (\pi, \delta)$, field of
fractions $K$ and  residue field 
$\kappa$. Suppose that $\kappa$ is a finite field of characteristic
not equal to 2.  Let $q_1 = <a_1,  \cdots , a_6>$ and $q_2 = <b_1,
b_2, b_3>$ be    quadratic forms  over $K$ with $a_i, b_j  \in
A$ square free and $q_1$ anisotropic. Suppose that
 the only primes in $A$ which divide
$a_1\cdots a_6b_1b_2 b_3$ are at most $\pi$ and  $\delta$.
Suppose that the dimensions of the first residues  (cf.
\cite{Sc}, Chapter 6,  \S 2, 2.5, page 209) of $q_1 \perp -q_2$ 
at $\pi$ and $\delta$ are  at least 5 and the dimensions of the first
residues of  $q_1$ and $q_2$ at $\pi$ and $\delta$ are at least 1.  
Then there exists  $x \in A^6$ and $y \in A^3$  such that $q_1(x) = q_2(y)$ 
  is either  a unit in $A$  or  of the form $
u\pi  + v\delta$   for  some units $u, v \in A^*$.  
\end{lemma}

\begin{proof}   Since the only prime divisors of $a_1 \cdots a_6b_1b_2b_3$ are at 
most  $\pi$ and $\delta$,  we write    $q_1    \simeq
\phi_0 \perp \phi_1\pi \perp \phi_2\delta \perp \phi_3\pi\delta$ and
$q_2 \simeq \psi_0 \perp \psi_1 \pi \perp \psi_2 \delta \perp
\psi_3\pi\delta$, for some  quadratic forms $\phi_i$ and
$\psi_i$ which are regular on $A$.   Since  $q_1$ is
anisotropic,  $\phi_i$'s are anisotropic  over the residue field $\kappa$ and hence of dimension 
at most  2,  $\kappa$ being a finite field.  Since the first
residues of  $q_1 \perp -q_2$ at $\pi$ and $\delta$ are at least  5
dimensional, the dimensions of $\phi_0 \perp -\psi_0 \perp \phi_2
\delta \perp -\psi_2 \delta 
$ and $\phi_0 \perp -\psi_0 \perp \phi_1\pi \perp
-\psi_1\pi$ are at least 5.  In particular either the
dimension of $\phi_0 \perp -\psi_0$ is at least 3 or the dimensions of $\phi_1
\perp -\psi_1$ and $\phi_2 \perp -\psi_2$ are both at least 3.

Suppose that the dimension of $\phi_0 \perp -\psi_0$ is at least
3.   Since the dimension of $\phi_0$ is at most 2, the dimension of
$\psi_0$ is at least 1.  If   the dimension of $\phi_0$ is at
least 1, then $\phi_0$ and $\psi_0$  represents a common unit over $A$,
 which is a common unit value of $q_1$ and $q_2$. 
If the dimension of $\phi_0$ is zero,  then the dimension of
$\psi_0$ is 3.  Then,   $q_2 = \psi_0$ is isotropic over $A$
and represents every element of $K$.  In this case the
dimensions of $\phi_1$, $\phi_2$ are equal to 2 and they 
represent all units in $A$. 
In particular, $q_1$ represents $u\pi + v \delta$ for  any  pair of units $u, v
\in A$,   which is also a value of $q_2$.

Suppose that the dimension of $\phi_0 \perp -\psi_0$ is at most 2. Then
the dimensions of $\phi_1 \perp -\psi_1$ and $\phi_2 \perp -\psi_2$ are
at least 3 and hence isotropic.    Since the first residues of $q_2$ at $\pi$
and $\delta$ are of dimensions at least 1, the dimensions of $\psi_1$
and $\psi_2$ are at most 2.   Hence $\phi_1$, $\phi_2,
\psi_1, \psi_2$ are of dimension at least 1 and  $\phi_1 \perp -\psi_1$, 
$\phi_2 \perp -\psi_2$ are isotropic. Hence  $\phi_1 $ and 
$\psi_1$ represent a common unit $u$ and  $\phi_2$ and $\psi_2$ represent a
common unit $v$ over $A$.   In particular  $q_1$ and $q_2$
 represent $u\pi + v \delta$ over $A$.  \end{proof}

\begin{prop}
\label{app-prop}
(cf. \cite{PS2}, Proposition 4.5)  Let $K$ be a
function field of a curve 
over a $p$-adic field $k$. Suppose  that $p \neq 2$. 
Let $\zeta = (f) \cdot (a) \cdot (b) \in H^3(K, \mu_2)$ and   $c_1,
c_2, c_3 \in K^*$.  Then   there exists  $g \in K^*$   which is  a value 
of the
quadratic form $<c_1, c_2, c_3 >$ such that $\zeta = (f) \cdot (g) \cdot
(h)$ for some $h \in K^*$. 
\end{prop}

\begin{proof}  Let ${\OO}$ be the ring of integers in
$k$. Then there exists a regular excellent  proper integral two
dimensional scheme ${\XX}$ over ${\OO}$  with function field $K$
and the support of  $f, a, b, c_1, c_2, c_3$ on ${\XX}$ is a
 union of regular curves
 $C_1, \cdots,
C_n$ with normal crossings. 
  Let ${\PP}$ be a  finite set of closed
points of ${\XX}$ containing   $C_i \cap C_j$ for all $i \neq
j$ and at least one closed point from each of $C_i$.  
 Let $B$ be the semi-local ring
at ${\PP}$. Then $B$ is a regular semi-local ring and hence a
unique factorisation domain.  For $1 \leq i \leq n$, let $\pi_i \in B$
be a prime defining $C_i$ on $B$.  
For each $P \in \PP$, let $B_P$ denote the local ring
at $P$. If $P \in C_i \cap C_j$ for some $i \neq j$, then
the maximal ideal $m_P$ of $B_P$ is generated by $\pi_i$
and $\pi_j$.   Without loss of generality,  we
assume that $f, a, b, c_1, c_2, c_3 
\in B$ and are square free. Since the support of each of $f$, $a$, $b$, $c_1$, $c_2$,
$c_3$  on $B$ is the union of $C_i$, the only primes in $B$ which divide any
of $f, a, b, c_1, c_2, c_3$ are  among  $\pi_i$, $1 \leq i \leq n$.  

For $1 \leq i \leq n$,  choose $\epsilon_i \in \{ 0 , 1 \}$ as follows: If
$\pi_i$ does not divide $c_1c_2c_3$,  then $\epsilon_i = 0$.  Suppose 
that $\pi_i$ divides $c_1c_2c_3$. If $\pi_i$ does not divide $f$, then
$\epsilon_i = 1$.  If $\pi_i$ divides $f$ and $\pi_i^2$ does not
divide $c_1c_2c_3$, then $\epsilon_i = 0$.  If $\pi_i$ divides $f$ and
$\pi_i^2$  divides  $c_1c_2c_3$, then $\epsilon_i = 1$.  Let $\theta =
\prod_{i= 1}^n \pi_i^{\epsilon_i}$. 

 Let $q_1 = <1, -f><-a, -b, ab>$, $q_2 = <c_1, c_2, c_3>$,  $q_1' =
 \theta q_1$ and $q_2' = \theta q_2$. 
 Let $C = C_i$ for some $i$. 
  Let $q'_{1C}$ and $q'_{2C}$
 be the first residues of $q'_{1C}$ and $q'_{2C}$ at $C$ respectively. 
 It is routine to check that $q'_{1C}$ and $q'_{2C}$ are  at least
 one dimensional and  $q'_{1C} \perp -q'_{2C}$  is   of dimension at least 5.
 Let $P \in C \cap \PP$ and $\hat{B}_P$ the completion of
 $B_P$ at its maximal ideal.
  Then, by (\ref{app-lemma2}),
 there exists $x_P \in \hat{B}^6_P$ and $y_P \in \hat{B}_P^3$
 such that $q'_1(x_P) = q'_2(y_P)$ is 
 either unit or $u_i\pi_i + u_j\pi_j$ for some units $u_i, u_j \in \hat{B}_P$.
  For each $i$, $ 1 \leq i \leq n$, let $x_{i_P} \in (\hat{B}_P/(\pi_i))^6$ and 
  $y_{i_P} \in (\hat{B}_P/(\pi_i))^3$ be the
 images of $x_P$ and $y_P$ respectively.  
 
 Suppose that   $\kappa(C_i)$ is a $p$-adic field.
 Then $P$ is the only closed point of $C_i$ and 
 $\kappa(C_i)$ is the field of fractions of $\hat{B}_P/(\pi_i)$. 
 Choose   $y_{i} =   x_{i_P}$.  
  Suppose that $\kappa(C_i)$
 is the function field of a curve over a finite field. 
 Let $q_{C_i} = \overline{q'}_1 - \overline{q'}_2$ be the 
 image of the quadratic form $q'_1 - q'_2$ in $B/(\pi_i)$.
 Note that $q_{C_i}$    could possible be singular and its non-singular part
 is $q'_{1C} - q'_{2C}$.  Since the first residue of  $q'_1 - q'_2$ at $C_i$ is at least 
 5 dimensional,  the non-singular part of $q_{C_i}$ is isotropic over $\kappa(C_i)$.  
  For each $P \in C_i \cap \PP$,   $(x_{i_P}, y_{i_P}) \in (\hat{B}_P/(\pi_i))^9$ 
  is an isotropic vector for $q_{C_i}$. By the weak approximation for (isotropic)
  quadrics,  there exist  vectors $x_i \in \kappa(C_i) ^6$ and $y_i 
  \in \kappa(C_i)^3$  
   such that $x_i$ and $y_i$ are  `sufficiently'
 close to $x_{i_P}$ and $y_{i_P}$ respectively for all $P \in C_i \cap \PP$.
 Since $x_{i_P} \in( \hat{B}_P/(\pi_i))^6$ for any $P \in C_i \cap \PP$ 
 and  is close to $x_i$, it follows that $x_i \in (B_P/(\pi_i))^6$
for all  $P \in C_i \cap \PP$ and hence $x_i \in (B/(\pi_i))^6$.
Similarly $y_i  \in (B/(\pi_i))^3$. Since 
$y_{i}$ is close to $y_{i_P}$ modulo $m_P$ for all $P \in C_i \cap \PP$,
$y_{i} =  y_{i_P}$ modulo $m_P$.  Since 
 $y_{i_P}$ is the image of $y_P$,  $y_{i} = y_{j}$ modulo $m_P$
 for  $P \in C_i \cap C_j$ for all $i \neq j$.  Thus, by (\cite{S3}, Proposition 0.3), 
 there exists a vector $y \in B^3$ such that $y = y_{i}$ modulo $\pi_i$
 for all $i$. 
 
Let $g' = q_2(y)$ and $g  = \theta g'$.
We claim that $(f) \cdot (a) \cdot (b) = (f) \cdot (g) \cdot (h)$
for some $h \in K^*$. 

Let  $\alpha = (f) \cdot (g)$.
Then, to prove our claim, it is enough to show that for 
each $z \in \XX$, there exists $h_z \in K_z^*$ such that 
$\zeta  - \alpha \cdot (h_z) \in H^3_{nr}(K_z/\hat{A}_z, \mu_l)$.

It follows from
 by (\cite{Sc}, Chapter 4, p.145, Theorem 1.7)
 that   over any field $(f) \cdot (a) \cdot (b) = (f) \cdot (g) \cdot (h)$
if and only if $g$ is a value of $q_1$. 
 
 Suppose that $z$ is a closed point of $\XX$. If $z \not\in C_i \cap C_j$
 for all $i \neq j$, then by (\ref{ramified-at-pi}), $\zeta$ is unramified 
 on $\hat{A}_z$ and hence $\zeta - \alpha \cdot (1) \in H^3_{nr}(K_z/
 \hat{A}_z, \mu_l)$.  Suppose that  $z \in    C_i \cap C_j$ for some $i \neq j$.
 By the choice of $y_{i}$ and $y_P$,  $g' \in B_P$ is either 
 a unit or $u_P\pi_i$ modulo $\pi_j$ or
 $v_P\pi_j$ modulo $\pi_i$ for some units $u_P, v_P$ at $P$.
 Since $q'_1 \perp <-g'>$ is isotropic over $K_{\pi_i}$
 and $K_{\pi_j}$, by (\ref{app-lemma1}), $q'_1 \perp <-g'>$
 is isotropic over $K_P$. Hence $g'$ is value of $q'_1$ over $K_P$.
 In particular, $g = \theta g'$ is a value of
 $q_1 = \theta q'_1$ over $K_P$ and  $(f) \cdot (a) \cdot (b) = (f) \cdot (g) \cdot (h_P)$
 for some $h_P \in K_P^*$.

Suppose that  $z$ is a codimension one point of $\XX$.
If $z$ is not the generic point of $C_i$ for some $i$, then
$\zeta$ is unramified at $z$ and hence 
$\zeta - \alpha \cdot (1) \in H^3_{nr}(K_z/\hat{A}_z, \mu_l)$.
Suppose that $z$ is the generic point of $C_i$ for some $i$.
By the choice of $y$, $g' = q'_2(y)$ is a value of 
$q'_1$ over $K_{C_i} = K_z$ and hence 
$g$ is a value of $q_1$.  Thus $\zeta = \alpha \cdot (h_z)$ for some
$h_z \in K_z^*$.  Hence, 
by (\ref{local-global-special}),  there exists $h \in K^*$ such that 
 $ (f) \cdot (a) \cdot (b)  - (f) \cdot (g) \cdot (h) \in H^3_{nr}(K/\XX, \mu_l)$.
 Since $H^3_{nr}(K/{\XX}, \mu_2) = 0$ (\cite{Ka}, 5.2), 
we have $\zeta = \alpha \cdot (h) = (f) \cdot (g) \cdot (h)$.
 \end{proof}

\begin{cor} 
\label{app-cor}
 (cf., \cite{PS2}, 4.6) Let $k$ be a $p$-adic field and $K$ be a
function field of a curve over $k$.  If $p \neq 2$, then $u(K) = 8$. 
\end{cor}

\begin{proof}    The proof follows as in (\cite{PS2}, 4.6) using
\ref{app-prop} and (\cite{PS1}, 4.4). \end{proof} 
 
\providecommand{\bysame}{\leavevmode\hbox to3em{\hrulefill}\thinspace}

\end{document}